\documentclass[a4paper, 11pt]{article}
\usepackage{amsmath,amssymb,esint,amscd,xspace,fancyhdr,color,authblk,srcltx,fontenc,bbm}
\usepackage{hyperref}
\usepackage{cite}

\newtheorem{theorem}{Theorem}[section]

\newtheorem{definition}[theorem]{Definition}
\newtheorem{lemma}[theorem]{Lemma}

\newtheorem{remark}[theorem]{Remark}

\newenvironment{proof}[1][Proof]{\textbf{#1.} }{\hfill\rule{0.5em}{0.5em}}
{\catcode`\@=11\global\let\AddToReset=\@addtoreset
\AddToReset{equation}{section}

\AddToReset{theorem}{section}

\title{\textsc{A large scaling property of level sets for degenerate $p$-Laplacian equations with logarithmic BMO matrix weights}}

\author{Thanh-Nhan Nguyen\thanks{Group of Analysis and Applied Mathematics, Department of Mathematics, Ho Chi Minh City University of Education, Ho Chi Minh City, Vietnam; \texttt{nhannt@hcmue.edu.vn}}, Minh-Phuong Tran\footnote{Corresponding author} \thanks{Applied Analysis Research Group, Faculty of Mathematics and Statistics, Ton Duc Thang University, Ho Chi Minh City, Vietnam; \texttt{tranminhphuong@tdtu.edu.vn}}}

\date{\today}

\begin{document}
 
\maketitle
\begin{abstract}
In this study, we deal with generalized regularity properties for solutions to $p$-Laplace equations with degenerate matrix weights. It has been already observed in previous interesting works~\cite{BDGN2022, BBDL2023} that gaining Calder\'on-Zygmund estimates for nonlinear equations with degenerate weights under the so-called $\log$-$\mathrm{BMO}$ condition and minimal regularity assumption on the boundary. In this paper, we also follow this direction and extend general gradient estimates for level sets of the gradient of solutions up to more subtle function spaces. In particular, we construct a covering of the super-level sets of the spatial gradient $|\nabla u|$ with respect to a large scaling parameter via fractional maximal operators.\\


\medskip

\end{abstract}

\maketitle


\section{Introduction}\label{sec:intro}

1.1. \textbf{Motivation and setting of the problems.}  
In this paper, we study the generalized regularity of weak solutions to the following non-homogeneous Dirichlet boundary value problem
\begin{align}\label{Eq-0}
\begin{cases}
-\mathrm{div}\left(\mathcal{L}_p(x,\nabla u)\right) &= \ -\mathrm{div}\left(\mathcal{L}_p(x,\mathbf{F})\right),  \quad \mbox{in} \ \Omega,\\
{\hspace{1cm}} u &= \ \mathsf{g}, \quad \mbox{on} \ \partial\Omega,
\end{cases}
\end{align}
where $\Omega \subset \mathbb{R}^n$ is a domain (=open subset) with non-smooth boundary $\partial\Omega$, for $n \ge 2$; $\mathbf{F}: \Omega \to \mathbb{R}^n$ is a given vector field, $\mathsf{g}: \overline{\Omega} \to \mathbb{R}$ is a measurable boundary datum. Further, we will motivate our interest in the specific case $\mathcal{L}_p(x,\zeta):= |\mathbb{P}(x)\zeta|^{p-2}\mathbb{P}^2(x)\zeta$, for $1<p<\infty$, that represents the elliptic operator driven by $p$-Laplacian involving the imposed matrix-valued weight $\mathbb{P}: \Omega \to \mathbb{R}^{n \times n}_{\mathrm{sym}^+}$ that is symmetric, positive definite and satisfies 
\begin{align}
\label{cond-PI}
|\mathbb{P}(x)|\cdot|\mathbb{P}^{-1}(x)| \le \Lambda, \quad x \in \Omega,
\end{align}
for some $\Lambda\ge 1$ and $|\cdot|$ denotes the default matrix norm induced by the Euclidean vector norm. Let us further define a scalar weight $\omega: \Omega \to [0,\infty)$ as follows
\begin{align}\label{def-omega}
\omega(x) = |\mathbb{P}(x)|, \quad x \in \Omega,
\end{align}
and the condition~\eqref{cond-PI} can be rewritten by
\begin{align}\label{Pi-2}
\Lambda^{-1} \omega(x) \mathrm{Id}_n \le \mathbb{P}(x) \le \omega(x) \mathrm{Id}_n,  \ \mbox{ or } \
\Lambda^{-1} \omega(x) |\zeta| \le |\mathbb{P}(x)\zeta| \le \omega(x) |\zeta|, 
\end{align}
for all $x \in \Omega$ and $\zeta \in \mathbb{R}^n$. Here, $\mathrm{Id}_n$ denotes the $n \times n$ identity matrix. Suppose that $\omega^p$ belongs to the class of Muckenhoupt weights $\mathrm{A}_p$, it is possible to consider the corresponding weighted Lebesgue space in the multiplicative sense $L^p_\omega(\Omega):=L^p(\Omega,\omega^pdx)$ and the corresponding Sobolev space $W^{1,p}_{\omega}(\Omega)$ (Muckenhoupt weights behave in a multiplicative form, see Section~\ref{sec:pre} for detailed definitions). Considering $\mathbf{F} \in L^{p}_{\omega}(\Omega)$ and $\mathsf{g} \in W^{1,p}_{\omega}(\overline{\Omega})$, we say that a weak solution to such problem~\eqref{Eq-0} is a map $u \in \mathsf{g} + W^{1,p}_{0,\omega}(\Omega)$ satisfying the weak formulation
\begin{align}\label{var-P0}
& \int_{\Omega} |\mathbb{P}(x)\nabla u|^{p-2}\mathbb{P}(x)\nabla u \cdot  \mathbb{P}(x)\nabla\varphi dx  =  \int_{\Omega} |\mathbb{P}(x)\mathbf{F}|^{p-2}\mathbb{P}(x)\mathbf{F} \cdot \mathbb{P}(x)\nabla \varphi dx,
\end{align}
for all $\varphi \in W^{1,p}_{0,\omega}(\Omega)$.

Equation~\eqref{Eq-0} appears naturally in different contexts as well as in variational models for many problems from mathematical physics. When $\mathbb{P}\equiv\mathrm{Id}_n$, the equation~\eqref{Eq-0} reads a non-homogeneous $p$-harmonic function, and to our knowledge, this type of operator appears a lot in physics, especially in the radiation of heat, glaciology, rheology, plastic molding, etc. Otherwise, in the case when $p=n$ and $\mathbf{F}=0$, the equation plays a crucial role in the theory of quasiconformal mappings, an important subject in complex analysis, as well as in physics and engineering. Note that by rewriting~\eqref{Eq-0} as
\begin{align*}
-\mathrm{div}\left((\mathbb{P}^2(x) \nabla u \cdot \nabla u)^{\frac{p-2}{2}}\mathbb{P}^2(x) \nabla u \right) = -\mathrm{div} \left( |\mathbb{P}(x)\mathbf{F}|^{p-2}\mathbb{P}^2(x)\mathbf{F}\right),
\end{align*}
it emphasizes that the latter is the Euler-Lagrange equation of minimizers of the functional
\begin{align*}
\mathcal{F}(w):=\int_\Omega{|\mathbb{P}(x) \nabla w|^p dx} - p \int_\Omega{|\mathbb{P}(x) \mathbf{F}|^{p-2}\mathbb{P}(x)\mathbf{F} \cdot (\mathbb{P}(x) \nabla w)dx}.
\end{align*}

During the last few years, there have been extensive mathematical investigations of both solvability and regularity theory for various classes of problems whose nonlinearity is connected with a matrix weight. For instance, a lot of authors have been executing their research to analyze the linear case when $\mathcal{L}_p(x,\nabla u) \equiv \mathbb{A}(x)\nabla u$, where $\mathbb{A}: \Omega \rightarrow \mathbb{R}^{n \times n}_{\text{sym}}$ is a uniformly elliptic weight, i.e.
\begin{align}
\label{eq:uniform}
\lambda_1 |\zeta|^2 \le \langle \mathbb{A}(x)\zeta,\zeta \rangle \le \lambda_2 |\zeta|^2, \quad \forall x \in \Omega, \quad \forall \zeta \in \mathbb{R}^n.
\end{align}
In particular, we refer the reader to~\cite{Meyer1963, KZ1999, AQ2002, FFZ2013, BW2004, Byun2005} for local and global regularity results of this standard model. Otherwise, concerning the case when $\mathbb{A}$ is uniformly elliptic with degenerate weight, that is,
\begin{align}
\label{eq:degenerate}
\Lambda^{-1} \mu(x)|\zeta|^2 \le \langle \mathbb{A}(x) \zeta, \zeta \rangle \le \Lambda \mu(x)|\zeta|^2, \quad \forall x \in \Omega, \quad \forall \zeta \in \mathbb{R}^n,
\end{align}
for some non-negative weight function $\mu$. A lot of attention has been devoted lately to this degenerate elliptic class, and the question of optimal regularity properties has attracted and been studied by many authors in~\cite{FKS1982, CMP2018, CMP2019, BDGN2022, BBDL2023} under various types of assumptions for $\mathbb{A}$ and $\mu$. To be more specific, concerning the study of the quasilinear elliptic equations of the kind~\eqref{Eq-0} with degenerate ellipticity condition~\eqref{eq:degenerate}, in~\cite{FKS1982} authors proved that when $\mu$ belongs to a Muckenhoupt class $\mathrm{A}_2$ and the data $\mathbf{F}$ is nice enough, weak solution $u \in C^{0,\beta}$ for some $\beta>0$. Later, several authors in these last years have extended gradient regularity in (weighted) Lebesgue spaces: Cao \emph{et al.} in~\cite{CMP2018} concluded the local gradient estimates that $\|\nabla u\|_{L^q(\mu dx)} \le C \|\mathbf{F}\|_{L^q(\mu dx)}$, for every $q>1$ when $\mu \in \mathrm{A}_2$ and $\mathbb{A}$ has small $\mathrm{BMO}_\mu^2$ norm; then the validity of global estimates was obtained by Phan in~\cite{Phan2020} and it is naturally extended to the vectorial case in~\cite{CMP2019}. Recently, an interesting new type of local gradient regularity was successfully presented by Balci \emph{et al.} in~\cite{BDGN2022}, stated that $\|\nabla u\|_{L^q(\omega^q dx)} \le C \|\mathbf{F}\|_{L^q(\omega^q dx)}$, for every $q \in (1,\infty)$, where a new small BMO assumption is imposed on $\mathrm{log}\mathbb{A}$ instead of the small $\mathrm{BMO}_\mu$ norm conditions presented in~\cite{CMP2018} (recall that $\mu = \Lambda^{-1}\omega^2$). To proceed further in the investigation, the introduced $\log$-$\mathrm{BMO}$ condition on the weight is also aimed to address gradient bounds for weighted $p$-Laplacian equation~\eqref{Eq-0}, where $\mathcal{L}_p(x,\zeta)= |\mathbb{P}(x)\zeta|^{p-2}\mathbb{P}^2(x)\zeta$.  Motivated by this study, some new global results are allowed to extend under an appropriate additional assumption on the domain $\Omega$, through the works in~\cite{BBDL2023, BC2024}. 

Following a recent trend of such interesting results in the literature, regarding the $p$-Laplacian equation with degenerate weights~\eqref{Eq-0}, we continue these works to investigate the gradient regularity associated with weak solutions in some generalized function settings. Specifically, inspired by the recent results concerning local and global Calder\'on-Zygmund estimates for nonlinear elliptic equations with degenerate weights (under the smallness $\log$-$\mathrm{BMO}$ condition) in~\cite{BDGN2022, BBDL2023}, by the covering argument of super level sets, we prove the following global implication via the presence of fractional maximal operators (see Definition~\ref{def:Malpha}):
\begin{align*}
\mathbf{M}_\alpha \left(|\mathbb{P} \mathbf{G}|^p \right) \in \mathbb{S}(\Omega) \Rightarrow \mathbf{M}_\alpha \left(|\mathbb{P} \nabla u|^p \right) \in \mathbb{S}(\Omega),
\end{align*}
under the minimal hypotheses of matrix weight $\mathbb{P}$ and $\partial\Omega$, where $\mathbf{G}:=\mathbf{F} + \nabla\mathsf{g}$ provides the information of given data of the problem, $\mathbf{M}_\alpha$ denotes the fractional maximal functions for $0<\alpha<n$, and moreover, many of the properties of classic Lebesgue spaces will be inherited by generalized ones described by $\mathbb{S}$ (see Section~\ref{sec:pre}). Our aim is, in particular, we deal with some instances of rearrangement-invariant quasi-normed spaces such as weighted generalized Lorentz, or generalized Morrey spaces.

Let us discuss some related issues that we believe might be meaningful to investigate regularity properties of solutions in terms of $\mathbf{M}_\alpha$. To the best of our knowledge, as shown in~\cite{Gilbarg1983, H1995, Hed1972}, the operator $\mathbf{M}_\alpha$ has a close relation with Riesz potential $\mathbf{I}_\alpha$ and from that point, via the fractional maximal operators of gradient of a function, it allows us to control the information of both size and oscillation of that function in the Lebesgue sense. We refer the reader in particular to the discussion in~\cite{KS2003, NPV2011} and the references therein concerning the detailed connection between $\mathbf{M}_\alpha$ and the so-called \emph{fractional derivatives} in fractional Sobolev spaces.

1.2. \textbf{Notation and main assumptions.} Before stating our results more precisely, let us specify some notation and assumptions. In the whole paper, for any $\rho>0$ and $\zeta \in \mathbb{R}^n$, we denote by $B(\zeta,\rho)$ the open ball in $\mathbb{R}^n$ of radius $\rho$ and center $\zeta$. We additionally write
$$\lambda B(\zeta,\rho) := B(\zeta, \lambda \rho), \mbox{ and } \lambda \Omega(\zeta,\rho) := \Omega \cap B(\zeta, \lambda \rho), \ \mbox{ for } \lambda \in \mathbb{R}^+.$$ 
In the following, we shall adopt the customary convention of denoting a constant by $C$, whose value is larger than one. Through estimates, $C$ may change the value from one line to another, and the dependencies of $C$ on prescribed parameters, if needed, will be kept between parentheses, sometimes will be properly emphasized at the end of the statements, for the sake of readability. For simplicity, we shall write $\mathsf{f} \in \mathcal{M}eas(\Omega,\mathbb{R})$ to indicate the Lebesgue measurable function $\mathsf{f}: \Omega \to \mathbb{R}$; $D_0=\mathrm{diam}(\Omega):=\sup_{y,z \in \Omega}{|y-z|}$ to employ the diameter of $\Omega$; and $\{|\mathsf{f}|>\sigma\} := \{\zeta \in \Omega: \, |\mathsf{f}(\zeta)|>\sigma\}$ in the arguments. Moreover, for a given measurable open set $\mathcal{O}$ of $\mathbb{R}^n$, we shall write $|\mathcal{O}|$ to mean the Lebesgue measure of $\mathcal{O}$ and $\overline{\mathsf{f}}_{\mathcal{O}} := \fint_{\mathcal{O}} \mathsf{f}(\zeta) d\zeta = \frac{1}{|\mathcal{O}|} \int_{\mathcal{O}} \mathsf{f}(\zeta) d\zeta$ the average value of every measurable function $\mathsf{f} \in L^1(\mathcal{O})$.

In what follows, we stress that we will try to use the notation $\mathbb{R}^{n \times n}_{\mathrm{sym}}$ to mean the set of all $n \times n$  symmetric matrices in $\mathbb{R}$; and $\mathbb{R}^{n \times n}_{\mathrm{sym}^+}$ denotes the subset of symmetric and positive definite matrices. In addition, for each matrix $M \in \mathbb{R}^{n \times n}_{\mathrm{sym}}$, we shall denote by $|M| = \sup_{|\zeta|\le 1} |M\zeta|$ to represent its spectral norm. As far as we are concerned, it is clear to define the mapping $\exp: \mathbb{R}^{n \times n}_{\mathrm{sym}} \to \mathbb{R}^{n \times n}_{\mathrm{sym}^+}$ and its inverse $\log: \mathbb{R}^{n \times n}_{\mathrm{sym}^+} \to \mathbb{R}^{n \times n}_{\mathrm{sym}}$ by using Taylor's theorem. 

Next, we shall shed some light on the main assumptions required for the given data of our problem. \\
\textbf{Assumption~$\mathbf{A_1}$} \textit{(Small $\log$-$\mathrm{BMO}$ condition)} Let $\mathbb{P}: \Omega \to \mathbb{R}^{n \times n}_{\mathrm{sym}^+}$ be a degenerate elliptic matrix-valued weight with uniformly bounded condition number as in~\eqref{cond-PI}. For given $R>0$, we define the $\log$-$\mathrm{BMO}$ semi-norm of a matrix weight $\mathbb{P}$  as follows
\begin{align}\label{BMO-norm}
|\log \mathbb{P}|_{\mathrm{BMO}(\mathbb{R}^n)} :=\sup_{z \in \mathbb{R}^n}{\sup_{0<r \le R} \fint_{B(z,r)} |\log \mathbb{P}(x) - \langle \mathbb{P}\rangle_{B(z,r)}^{\log}| dx}, 
\end{align}
where $\langle \mathbb{P}\rangle_B^{\log}$ is the logarithm average of $\mathbb{P}$ over the ball $B$, given as
\begin{align}\label{log-PI}
\langle \mathbb{P}\rangle_B^{\log} := \exp \left(\fint_B \log \mathbb{P}(x) dx\right).
\end{align}
Here, for the sake of brevity, in case the ball $B$ covers $\Omega$, we omit $B$ by just writing $|\log \mathbb{P}|_{\mathrm{BMO}}$ when no confusion arises. For given $\kappa>0$ and $R>0$, we say that the weight $\mathbb{P}$ satisfies the \emph{small $(\kappa,R)$-$\log$-$\mathrm{BMO}$ condition}, or equivalently, the $\log\mathbb{P}$ is \emph{$(\kappa,R)$-small-$\log$-$\mathrm{BMO}$} if $|\log \mathbb{P}|_{\mathrm{BMO}} \le \kappa$. As far as we know, regarding Assumption~$\mathbf{A_1}$, the idea of smallness-BMO condition on the logarithm of the weight was first mentioned by Balci \emph{et al.} in~\cite{BDGN2022} and here, our strategy is based on making use of this assumption to deal with the upper level set for gradient of solutions to~\eqref{var-P0}. \\[5pt]
\textbf{Assumption~$\mathbf{A_2}$} \textit{(The $(\kappa,r_0)$-Lipschitz condition for the boundary $\partial\Omega$)} Let $\kappa \in \left[0,\frac{1}{2n}\right]$ and $r_0>0$ be given. We say that $\Omega$ satisfies the $(\kappa,r_0)$-Lipschitz condition if and only if for every $x_0 \in \partial \Omega$, there exists a coordinate system $\{y_1, y_2,...,y_n\}$ such that $x_0$ is the origin in this system and a Lipschitz function $\Upsilon: \mathbb{R}^{n-1} \to \mathbb{R}$ satisfying $\|\nabla \Upsilon\|_{L^{\infty}} \le \kappa$ and
\begin{align*}
B(x_0,r_0) \cap \Omega = \left\{(y_1, y_2,..., y_n) \in B(x_0,r_0): \ y_n > \Upsilon(y_1,y_2,...,y_{n-1}) \right\}.
\end{align*}
{\bf Remark. } As shown in~\cite{BBDL2023}, the Lipschitz condition imposed on domain $\Omega$ is sharp. We also refer the interested reader to~\cite[Example 4.1]{BBDL2023}, in which the authors carefully provided a two-dimensional example to show that the Assumption~$\mathrm{A}_2$ on domain $\Omega$ is optimal when concluding the local regularity estimates imply the global ones. Therefore, as one could expect, it sufficiently allows us to obtain global results concerning problems with degenerate matrix weights.\\

1.3. \textbf{Statements of main results.} With these standing assumptions at hand, we are now in the position to state our main results. Besides, for the sake of brevity, the structural data of the problem will be deliberately not repeated in our statements and proofs in the paper. We use the abbreviation $\mathtt{data}$ to indicate the set of specified constants as follows
\begin{align*}
\mathtt{data} \equiv \texttt{data}\left(n,p,D_0,\alpha,[\mu]_{\mathrm{A}_{\infty}},\Lambda,r_0\right).
\end{align*}

The first theorem plays a key role in our study, which states the large-scaling level-set inequality involving the weighted fractional maximal distribution functions $d^\mu_\alpha$. For the readers' convenience, we also highlight this function here and the reader is forwarded to Definition~\ref{def:WFMD} of Section~\ref{sec:pre} for its detailed definition. 
\begin{align*}
d^\mu_\alpha(\mathsf{f},\lambda) = \int_{\{\mathbf{M}_{\alpha}\mathsf{f}>\lambda\}} \mu(x) dx, \quad \lambda > 0,
\end{align*}
for each $\alpha \in [0,n]$ and $\mathsf{f} \in L^1_{\mathrm{loc}}(\mathbb{R}^n)$. It emphasizes that the idea of using this term stemmed from our previous work~\cite{PNJFA} when we wanted to discuss the unified approach to the regularity via $\mathbf{M}_\alpha$. Moreover, this technique is based on the effective \emph{Harmonic free} method that lies at the heart of Acerbi-Mingione's work in~\cite{AM2007}, Byun-Wang's in~\cite{BW2004}.  
\begin{theorem}\label{theo:A}
Let $\mathbb{P}$ be a matrix weight satisfying~\eqref{cond-PI}; $\mathbf{F} \in L^{p}_{\omega}(\Omega)$ and $\mathsf{g} \in W^{1,p}_{\omega}(\overline{\Omega})$ for $p \in (1,\infty)$. Assume further that $u \in \mathsf{g} + W^{1,p}_{0,\omega}(\Omega)$ is a weak solution to~\eqref{Eq-0}, a Muckenhoupt weight $\mu \in \mathrm{A}_{\infty}$ and $\alpha \in [0,n)$. Then, for every $\varepsilon>0$ small enough  and $\theta>0$, one can find some positive constants 
$$\gamma = \gamma(\mathtt{data},\theta), \ \mbox{ and } \ \kappa = \kappa(\mathtt{data},\theta,\varepsilon)$$ 
such that if $\log \mathbb{P}$ is $(\kappa,r_0)$-small-$\log$-$\mathrm{BMO}$ and $\Omega$ satisfies $(\kappa,r_0)$-Lipschitz condition for some $r_0>0$, then the following estimate of the type
\begin{align}\label{dist-ineq}
 & d^{\mu}_{\alpha}(|\mathbb{P}(x)\nabla u|^{p};\varepsilon^{-\theta}\lambda) \le C \varepsilon d^{\mu}_{\alpha}(|\mathbb{P}(x)\nabla u|^{p};\lambda) + d^{\mu}_{\alpha}(|\mathbb{P}(x)\mathbf{G}|^{p};\varepsilon^{\gamma}\lambda)
\end{align}
holds for all $\lambda>0$. Here, the positive constant $C$ depends on $\mathtt{data}, \gamma$.
\end{theorem}

The main idea of the proof goes back to previous approaches mentioned above, on the one hand, allow us to prove suitable comparison estimates in our problems with degenerate weight (both local interior and up-to-boundary estimates), and on the other hand, allow us to combine various nontrivial covering techniques. The key point to the proof of Theorem~\ref{theo:A} essentially combines the techniques introduced in~\cite{PNJFA, AM2007, BW2004} with some novel insights that allow us to analyze the level sets of the fractional maximal function of the spatial gradient $|\nabla u|$.

In the next theorem, we state regularity estimates for gradients of weak solutions to~\eqref{Eq-0} in various generalized function settings, that could be useful for several purposes, for example in assessing the convergence of some optimization algorithms for min/max problems;  the gradient norm provides crucial information about the direction and rate of change of the energy functional, which is invaluable for the optimization process towards convergence; or in machine learning and deep learning, it provides the information about how steep the function is at a given point in the space, etc. 

Let us stress the reader's attention to the fact that Theorem~\ref{theo:general} here provides regularity results in a general form for the sake of completeness and the convenience of reading. We shall separate our statements and proofs in each desired function space estimate, which are presented in Section~\ref{sec:proof}. More precisely, we deal with some instances of rearrangement-invariant quasi-normed spaces such as weighted generalized Lorentz, or generalized Morrey spaces, etc.

\begin{theorem}\label{theo:general}
Let $\mathbb{P}$ be a matrix weight satisfying~\eqref{cond-PI}; $\mathbf{F} \in L^{p}_{\omega}(\Omega)$ and $\mathsf{g} \in W^{1,p}_{\omega}(\overline{\Omega})$ for $p \in (1,\infty)$. Assume further that $u \in \mathsf{g} + W^{1,p}_{0,\omega}(\Omega)$ is a weak solution to~\eqref{Eq-0}. Then, for any $\alpha \in [0,n)$, there exists $\kappa = \kappa(\mathtt{data},\mathtt{indices})>0$ such that 
\begin{align}\label{est:B}
\|\mathbf{M}_{\alpha}(|\mathbb{P}\nabla u|^{p})\|_{\mathbb{S}(\Omega)}  \le C\left\|\mathbf{M}_{\alpha}(|\mathbb{P}\mathbf{G}|^{p})\right\|_{\mathbb{S}(\Omega)},
\end{align}
 if $\log \mathbb{P}$ satisfies $(\kappa,r_0)$-small-$\log$-$\mathrm{BMO}$ condition and $\Omega$ is $(\kappa,r_0)$-Lipschitz domain for some $r_0>0$. Here, the simplified notation $\mathbb{S}$ employs relevant generalized function spaces with prescribed $\mathtt{indices}$; and constant $C$ depends on $\mathtt{data, indices}$. \end{theorem}

1.4. \textbf{Organization of the paper.} The introductory section is closed by highlighting the organization of the paper. First, some basic definitions and preliminary tools on matrix-valued weights, logarithms, and Muckenhoupt weights will be reviewed in the next section, Section~\ref{sec:pre}. Section~\ref{sec:lemma} consists of some preliminary lemmas that treat the comparison estimates for solutions in the interior and near the boundary points. Finally, the proofs of the main results are given in Section~\ref{sec:proof}. 

\section{Standard definitions and basic properties}
\label{sec:pre}
This section recalls several real analysis definitions and tools that are needed in the paper. We also focus our attention in the \emph{new} matrix weights and logarithms discussed in~\cite{BDGN2022, BBDL2023}. We first recall here the definition of $\mathrm{A}_p$-Muckenhoupt class of weights. 

\begin{definition}
Let $1 \le q \le \infty$ and a locally integrable function $\mu: \mathbb{R}^n \to [0,\infty)$. We say that this function belongs to the class of Muckenhoupt weights, i.e. $\mu \in \mathrm{A}_{q}$ for $q \in (1,\infty)$, if and only if the term $[\mu]_{\mathrm{A}_{q}}$ is finite, where
\begin{align*}
[\mu]_{\mathrm{A}_{q}} := \displaystyle \sup_{y \in \mathbb{R}^n, \varrho>0} \left(\fint_{B(y,\varrho)} [\mu(\zeta)]^{-\frac{1}{q-1}} d\zeta\right)^{q-1} \left(\fint_{B(y,\varrho)} \mu(\zeta) d\zeta\right).
\end{align*}
We say that $\mu \in \mathrm{A}_{1}$ if and only if there exists a constant $C$ such that
\begin{align*}
\sup_{\varrho>0} \fint_{B(x,\varrho)} \mu(\zeta) d\zeta \le C \mu(x), \mbox{ for a.e. } x \in \mathbb{R}^n.
\end{align*}
In this case, we define by $[\mu]_{\mathrm{A}_{1}}$ the smallest value of $C$ for which this inequality holds. On the other hand, we say that $\mu \in \mathrm{A}_{\infty}$ if and only if there exist  $C>0$ and $\nu$ such that the following estimate holds
\begin{align*}
\mu(\mathcal{O}) \le C \left({|\mathcal{O}|}/{|B|}\right)^\nu \mu(B),
\end{align*}
for all measurable subset $\mathcal{O}$ of any ball $B$. It also remarks that if $\mu \in L^1_{\mathrm{loc}}(\mathbb{R}^n;\mathbb{R}^+)$, we shall denote
\begin{align*}
\mu(\mathcal{O}) :=\int_{\mathcal{O}}\mu(x)dx,
\end{align*}
for some measurable set $\mathcal{O} \subset \mathbb{R}^n$. Furthermore, when $\mu \equiv 1$, it is known that $\mu(\mathcal{O}) \equiv |\mathcal{O}|$.
\end{definition}

\begin{remark}
It is worth mentioning here that $\mathrm{A}_{1} \subset \mathrm{A}_{q} \subset \mathrm{A}_{\infty}$ for any $1<q<\infty$ and $\mathrm{A}_{\infty}:= \cup_{q \ge 1} \mathrm{A}_{q}$. Moreover, $\mu \in \mathrm{A}_{\infty}$ if and only if there exist $c_1, c_2$ and $\nu_1,\nu_2>0$ such that
\begin{align}\label{Muck-w}
c_1 \left[{|\mathcal{O}|}/{|B|}\right]^{\nu_1} \mu(B)  \le \mu(\mathcal{O}) \le c_2 \left[{|\mathcal{O}|}/{|B|}\right]^{\nu_2} \mu(B). 
\end{align}
In this case, we shall write $[\mu]_{\mathrm{A}_\infty} = (c_1,c_2,\nu_1,\nu_2)$ for the sake of brevity.
\end{remark}

\begin{remark}
For any ball $B \subset \mathbb{R}^n$, we have another interesting property: $\mu^q \in \mathrm{A}_{q}$ if and only if two following inequalities hold true
\begin{align}\notag
\left(\fint_B [\mu(x)]^q dx \right)^{\frac{1}{q}} \le C_1 \langle \mu\rangle_B^{\log}, \quad \text{and} \quad \left(\fint_B [\mu(x)]^{-q'} dx \right)^{\frac{1}{q'}} \le C_2 \frac{1}{\langle \mu\rangle_B^{\log}},
\end{align}
where one further notices that
\begin{align*}
\langle 1/\mu\rangle_B^{\log} = \exp\left(-\fint_B \log \mu(x) dx\right) = \frac{1}{\langle \mu\rangle_B^{\log}}. 
\end{align*}
\end{remark}

As shown in~\cite{BBDL2023}, one concludes that if $|\log \mathbb{P}|_{\mathrm{BMO}}$ is small enough then $|\mathbb{P}|^q$ belongs to Muckenhoupt class $\mathrm{A}_q$. We restate it in the next lemma and the detailed proof can be found in~\cite{BDGN2022}.
\begin{lemma}\label{lem:Pi}
Let $1<q<\infty$, $\mathbb{P}: \Omega \to \mathbb{R}^{n \times n}_{\mathrm{sym}^+}$ be a matrix-valued weight function and $\mu = |\mathbb{P}|$. Then, $\mu^q \in \mathrm{A}_q$ when $\mathbb{P}$ satisfies the small $\log$-$\mathrm{BMO}$ condition. It means that, there exists a constant $\kappa = \kappa(q)>0$ such that if $|\log \mathbb{P}|_{\mathrm{BMO}} \le \kappa$ then $\mu^q \in \mathrm{A}_q$.
\end{lemma}

\begin{definition}[Weighted Lebesgue and Sobolev spaces] Let $1<q<\infty$ and $\mu$ be a given weight function, we define by
\begin{align*}
L^q_{\mu}(\Omega) \equiv L^q(\Omega,\mu^qdx) := \left\{\mathsf{f} \in \mathcal{M}eas(\Omega,\mathbb{R}) \mbox{ such that } \|\mathsf{f}\|_{L^q_{\mu}(\Omega)}< \infty\right\}
\end{align*}
the weighted Lebesgue space, where the term $\|\mathsf{f}\|_{L^q_{\mu}(\Omega)}$ is given by
$$\|\mathsf{f}\|_{L^q_{\mu}(\Omega)} := \left(\int_{\Omega} |\mathsf{f}(x)|^q [\mu(x)]^q dx\right)^{\frac{1}{q}}.$$ 
Moreover, the corresponding weighted Sobolev space is defined by
\begin{align*}
W^{1,q}_{\mu}(\Omega) = \left\{\mathsf{f} \in L^q_{\mu}(\Omega): \, |\nabla \mathsf{f}| \in L^q_{\mu}(\Omega) \right\}.
\end{align*}
The Sobolev space will be equipped with the following norm 
$$\|\mathsf{f}\|_{W^{1,q}_{\mu}(\Omega)} = \|\mathsf{f}\|_{L^q_{\mu}(\Omega)} + \|\nabla \mathsf{f}\|_{L^q_{\mu}(\Omega)}.$$ 
Here, we will also write $W^{1,q}_{0,\mu}(\Omega)$ to mean the closure of $C_0^{\infty}(\Omega)$ in $W^{1,q}_{\mu}(\Omega)$.
\end{definition}

As already mentioned, an important feature of our study is the role of fractional maximal operators $\mathbf{M}_\alpha$ in regularity estimates. We shall recall its definition here and note that the Hardy-Littlewood maximal function is just a specific form of such an operator. Moreover, Lemma~\ref{lem:bound-M-beta} presents the usual boundedness property, we refer to~\cite{PNJMAA} for detailed proof. 

\begin{definition}\label{def:Malpha}
The fractional maximal operator $\mathbf{M}_\alpha$ with $\alpha \in [0, n]$ is a measurable map defined on $L^1_{\mathrm{loc}}(\mathbb{R}^n)$ reads as
\begin{align}\label{eq:Malpha}
\mathbf{M}_\alpha \mathsf{f}(z) = \sup_{\varrho \in \mathbb{R}^+} {\varrho^{\alpha} \fint_{B(z,\varrho)}{|\mathsf{f}(\zeta)|d\zeta}},
\end{align}
for $z \in \mathbb{R}^n$, and $\mathsf{f} \in L^1_{\mathrm{loc}}(\mathbb{R}^n)$.
\end{definition} 

\begin{lemma}\label{lem:bound-M-beta}
For every $1 \le q < \infty$ and $\alpha \in \left[0,\frac{n}{q}\right)$, there exists $C(\alpha,q,n)>0$ such that
\begin{align*}
\sup_{\lambda \in \mathbb{R}^+} \lambda^{q} \left|\left\{\zeta \in \mathbb{R}^n: \, \mathbf{M}_{\alpha}\mathsf{f}(\zeta)>\lambda\right\}\right|^{1-\frac{\alpha q}{n}} \le C(\alpha,q,n) \int_{\mathbb{R}^n}|\mathsf{f}(\zeta)|^q d\zeta,
\end{align*}
for every $\mathsf{f} \in L^q(\mathbb{R}^n)$.
\end{lemma}

Let us recall the definition of distribution functions that correspond to a Muckenhoupt weight. The classical distribution function was introduced by Grafakos in~\cite{Grafakos}, and later, this concept can be modified and made use in our series of works in~\cite{NT20, PNJFA, MPT2018, PNJDE, NTT23, TN19-2, PNJMAA}. 

\begin{definition}[Distribution functions]
\label{def:WFMD}
Let $\mu \in \mathrm{A}_{\infty}$ and $\mathsf{f} \in \mathcal{M}eas(\Omega,\mathbb{R})$. The weighted distribution function of $\mathsf{f}$ over $\Omega$ is defined by
\begin{align}\label{def-Df}
d^{\mu}(\mathsf{f};\lambda) := \int_{\left\{|\mathsf{f}|> \lambda\right\}} \mu(x) dx, \quad \lambda \in \mathbb{R}^+.
\end{align}
Moreover, for $\mathsf{f} \in L^1_{\mathrm{loc}}(\Omega)$ and $\alpha \in [0,n]$, the weighted fractional distribution function $d^{\mu}_{\alpha}$ will be defined by
\begin{align}\label{def-Df-al}
d^{\mu}_{\alpha}(\mathsf{f};\lambda) := d^{\mu}(\mathbf{M}_{\alpha}\mathsf{f};\lambda), \quad \lambda \in \mathbb{R}^+.
\end{align}
\end{definition}

In this paper, we focus the study on the regularity estimates on some generalized function spaces that follow an interesting \emph{rearrangement invariant property}. 

\begin{definition}[Weighted Lorentz spaces] Let $\mu \in \mathrm{A}_{\infty}$ and two parameters $\mathfrak{q} \in (0,\infty)$, $0 < \mathfrak{s}\le \infty$. The weighted Lorentz space $L^{\mathfrak{q},\mathfrak{s}}_{\mu}(\Omega)$ is defined by
\begin{align*}
L^{\mathfrak{q},\mathfrak{s}}_{\mu}(\Omega) :=\left\{\mathsf{f} \in \mathcal{M}eas(\Omega,\mathbb{R}): \ \|\mathsf{f}\|_{L^{\mathfrak{q},\mathfrak{s}}_{\mu}(\Omega)}<\infty\right\},
\end{align*}
where $\|\mathsf{f}\|_{L^{\mathfrak{q},\mathfrak{s}}_{\mu}(\Omega)}$ is given by
\begin{align}\label{w-Lor-norm}
\|\mathsf{f}\|_{L^{\mathfrak{q},\mathfrak{s}}_{\mu}(\Omega)} :=  \begin{cases} \left( \displaystyle{\int_{\mathbb{R}^+} \mathfrak{q} \left[\lambda^{\mathfrak{q}} d^{\mu}(\mathsf{f};\lambda)\right]^{\frac{\mathfrak{s}}{\mathfrak{q}}} \frac{d\lambda}{\lambda}} \right)^{\frac{1}{\mathfrak{s}}} \ &\mbox{ if } \ \mathfrak{s} < \infty,\\  
\displaystyle{\sup_{\lambda \in \mathbb{R}^+} \, \left[\lambda^{\mathfrak{q}} d^{\mu}(\mathsf{f};\lambda)\right]^{\frac{1}{\mathfrak{q}}}} \ &\mbox{ if } \ \mathfrak{s} = \infty.\end{cases}
\end{align}
We remind that $d^{\mu}$ is the weighted distribution function defined as in~\eqref{def-Df}.
\end{definition}

\begin{definition}[Generalized Lorentz spaces involving two weights]
Let $\mu \in \mathrm{A}_{\infty}$ and $\nu \in L^1_{\mathrm{loc}}(\mathbb{R}^+;\mathbb{R}^+)$ be a new weight. We further introduce a non-decreasing function $\Sigma$ as following
\begin{align}\label{def:Gam}
\Sigma(\tau) = \int_0^{\tau} \nu(s) ds, \quad \tau \in [0,\infty).
\end{align}
For each pair $\mathfrak{q} \in (0,\infty)$, $0<\mathfrak{s} \le \infty$ and $\mathsf{f} \in \mathcal{M}eas(\Omega,\mathbb{R})$, we will denote
\begin{align}\label{GL-norm}
\|\mathsf{f}\|_{L^{\mathfrak{q},\mathfrak{s}}_{\mu,\nu}(\Omega)} := \begin{cases}  \left(\displaystyle{\int_{\mathbb{R}^+} \mathfrak{q} \left[\lambda^{\mathfrak{q}}\Sigma\left(d^{\mu}(\mathsf{f};\lambda)\right) \right]^{\frac{\mathfrak{s}}{\mathfrak{q}}} \frac{d\lambda}{\lambda}} \right)^{\frac{1}{\mathfrak{s}}} \ &\mbox{ if } \mathfrak{s}<\infty,\\ 
\displaystyle{\sup_{\lambda \in \mathbb{R}^+} \, \left[\lambda^{\mathfrak{q}}\Sigma\left(d^{\mu}(\mathsf{f};\lambda)\right)\right]^{\frac{1}{\mathfrak{q}}}} \ &\mbox{ if } \mathfrak{s}=\infty.\end{cases}
\end{align}
Then, the generalized weighted Lorentz spaces, often written by $L^{\mathfrak{q},\mathfrak{s}}_{\mu,\nu}(\Omega)$, is the set of all functions $\mathsf{f} \in \mathcal{M}eas(\Omega,\mathbb{R})$ such that $\|\mathsf{f}\|_{L^{\mathfrak{q},\mathfrak{s}}_{\mu,\nu}(\Omega)}<\infty$.
\end{definition}

\begin{definition}[Generalized $\psi$-Morrey spaces] Let $\psi \in \mathcal{M}eas(\Omega \times \mathbb{R}^+,\mathbb{R}^+)$ and $\mathfrak{q} \in (0, \infty)$. The generalized Morrey space $\mathrm{M}^{\mathfrak{q},\psi}(\Omega)$ is defined by
\begin{align*}
\mathrm{M}^{\mathfrak{q},\psi}(\Omega) := \left\{\mathsf{f} \in L^\mathfrak{q}(\Omega): \ \|\mathsf{f}\|_{\mathrm{M}^{\mathfrak{q},\psi}(\Omega)}<\infty\right\},
\end{align*}
where $\|\mathsf{f}\|_{\mathrm{M}^{\mathfrak{q},\psi}(\Omega)}$ is given as
\begin{align}\label{Morrey-norm}
\|\mathsf{f}\|_{\mathrm{M}^{\mathfrak{q},\psi}(\Omega)} :=  \sup_{z\in \Omega; \, 0<r<D_0} \left(\frac{1}{\psi(z,r)}\int_{\Omega(z,r)}|\mathsf{f}(\zeta)|^{\mathfrak{q}}d\zeta\right)^{1/\mathfrak{q}}.
\end{align}
\end{definition} 

\section{Technical lemmas}
\label{sec:lemma}

From this section onwards, the content takes on a more analytic flavor. We shall present and prove some auxiliary tools that play an important role in the rest of the paper. In addition, a series of comparison estimates to suitable reference problems in local interior and boundary of domain will be established. We first discuss on the existence of weak solutions to~\eqref{var-P0} and prove a very first global estimate for such solutions in the $L^p$-sense. 

\begin{lemma}\label{lem:global}
Let $\mathbb{P}$ be a matrix weight satisfying~\eqref{cond-PI} and $\omega$ be defined as in~\eqref{def-omega}. Assume that $\mathbf{F} \in L^{p}_{\omega}(\Omega)$ and $\mathsf{g} \in W^{1,p}_{\omega}(\overline{\Omega})$ with given $p>1$. Then, there exists a small constant $\kappa>0$ such that if $|\log \mathbb{P}|_{\mathrm{BMO}} \le \kappa$ and equation~\eqref{Eq-0} admits a weak solution $u \in \mathsf{g} + W^{1,p}_{0,\omega}(\Omega)$. Furthermore, there exists a constant $C = C(n,p,D_0,\Lambda)>0$ such that
\begin{align}\label{est-global}
\int_{\Omega} |\mathbb{P}(x) \nabla u|^p dx \le C \left( \int_{\Omega} |\mathbb{P}(x) \mathbf{F}|^p dx + \int_{\Omega} |\mathbb{P}(x) \nabla\mathsf{g}|^p dx\right).
\end{align}
\end{lemma}
\begin{proof}
Thanks to Lemma~\ref{lem:Pi}, there exists $\kappa>0$ such that if $|\log \mathbb{P}|_{\mathrm{BMO}} \le \kappa$ then $\omega^p \in \mathrm{A}_{p}$. Hence, the existence of a weak solution $u \in \mathsf{g} + W^{1,p}_{0,\omega}(\Omega)$ to~\eqref{Eq-0} is ensured for this small $\log$-$\mathrm{BMO}$ semi-norm of $\mathbb{P}$. The proof of~\eqref{est-global} is simple by testing $u-\mathsf{g}$ to~\eqref{var-P0} and applying Young's inequality.
\end{proof}

The following preliminary lemma is useful for our need later in comparison procedures. With regards to other related inequalities on the uniformly convex Orlicz functions, we also refer the reader to~\cite[Appendix B]{DFTW20}, where the authors carefully proved several notable results. 
\begin{lemma}\label{lem-Phi}
Let $p >1$ and two functions $\Psi: \mathbb{R}^+ \to \mathbb{R}^+$, $\mathbb{V}_{p}: \mathbb{R}^n \to \mathbb{R}^n$ be defined by
\begin{align}\label{def-Phi}
\Psi(t) := \frac{1}{p}t^p, \  t \in \mathbb{R}^+ \ \mbox{ and }  \ \mathbb{V}_{p}(\zeta) := |\zeta|^{\frac{p-2}{2}}\zeta, \ \zeta \in \mathbb{R}^n.
\end{align}
If $p \ge 2$ then there exists a constant $C = C(p)>0$ such that
\begin{align}\label{vphi-1}
\Psi(|\zeta_1-\zeta_2|) \le C |\mathbb{V}_{p}(\zeta_1) - \mathbb{V}_{p}(\zeta_2)|^2, \quad \forall \zeta_1, \zeta_2 \in \mathbb{R}^n.
\end{align}
Otherwise, if $1<p<2$, for every $\epsilon \in (0,1)$ there exists a constant $C>0$ such that
\begin{align}\label{vphi-2}
\Psi(|\zeta_1-\zeta_2|) \le \epsilon \Psi(|\zeta_1|) + C {\epsilon}^{1 -\frac{2}{p}} |\mathbb{V}_{p}(\zeta_1) - \mathbb{V}_{p}(\zeta_2)|^2, \quad \forall \zeta_1, \zeta_2 \in \mathbb{R}^n.
\end{align}
\end{lemma}
\begin{proof} 
Let us first recall the shifted $N$-function associated to $\Psi$ as below
\begin{align*}
\Psi_a(t) := \int_0^t \frac{s\Psi'(\max\{a;s\})}{\max\{a;s\}}ds, \quad a, t \in \mathbb{R}^+.
\end{align*}
For every $a, t \in \mathbb{R}^+$, by a simple computation, we can show that
\begin{align}\notag
\Psi_a(t) \simeq (\max\{a,t\})^{p-2}t^2.
\end{align}
That means there exist two positive constants $C_1, C_2>0$ such that
\begin{align}\label{simi}
C_1\Psi_a(t) \le (\max\{a,t\})^{p-2}t^2 \le C_2 \Psi_a(t), \quad \mbox{ for all } a,t \in \mathbb{R}^+.
\end{align}
The proof of~\eqref{vphi-1} is very simple for the first case $p \ge 2$. Indeed, by~\eqref{simi} one has
\begin{align}\label{phi-t}
\Psi(t) = \frac{1}{p}t^p = \frac{1}{p} t^{p-2}t^2 \le \frac{1}{p} \big(\max\{a,t\}\big)^{p-2} t^2 \le \frac{1}{p} C_2\Psi_a(t).
\end{align} 
Moreover, for all $\zeta_1, \zeta_2 \in \mathbb{R}^n$, it is well-known that 
$$|\mathbb{V}_{p}(\zeta_1) - \mathbb{V}_{p}(\zeta_2)|^2 \simeq \Psi_{|\zeta_2|}(|\zeta_1-\zeta_2|),$$ which means there exist $C_3, C_4>0$ such that
\begin{align}\label{Phi-PQ}
C_3\Psi_{|\zeta_2|}(|\zeta_1-\zeta_2|) \le |\mathbb{V}_{p}(\zeta_1) - \mathbb{V}_{p}(\zeta_2)|^2 \le C_4 \Psi_{|\zeta_2|}(|\zeta_1-\zeta_2|).
\end{align}
Therefore, one may obtain~\eqref{vphi-1} from~\eqref{phi-t} and~\eqref{Phi-PQ}. It is worth mentioning that all constants $C_1, C_2, C_3, C_4$ in~\eqref{simi} and~\eqref{Phi-PQ} only depend on $p$. \\

We now show~\eqref{vphi-2} for the remain case $1<p<2$. For all $\zeta_1, \zeta_2 \in \mathbb{R}^n$, we may use the decomposition
\begin{align}\notag
|\zeta_1-\zeta_2|^p & = \left(\big(|\zeta_1|+|\zeta_2|\big)^{p-2} |\zeta_1-\zeta_2|^2\right)^{\frac{p}{2}}  \big(|\zeta_1|+|\zeta_2|\big)^{\frac{(2-p)p}{2}} .
\end{align}
Combining two following basic inequalities 
\begin{align*}
\max\{|\zeta_1-\zeta_2|; |\zeta_2|\} \le |\zeta_1| + |\zeta_2| \mbox{ and } (|\zeta_1|+|\zeta_2|)^p \le 4^p\big(|\zeta_1-\zeta_2|^p + |\zeta_1|^p\big),
\end{align*}
one gets that
\begin{align}\notag
|\zeta_1-\zeta_2|^p & \le C \left(\big(\max\{|\zeta_1-\zeta_2|; |\zeta_2|\}\big)^{p-2} |\zeta_1-\zeta_2|^2\right)^{\frac{p}{2}} \big(|\zeta_1|^p+|\zeta_1-\zeta_2|^p\big)^{\frac{2-p}{2}}.
\end{align}
This inequality is equivalent to
\begin{align}\label{est-phi-PQ}
\Psi(|\zeta_1-\zeta_2|) \le C \big[\Psi(|\zeta_1|) + \Psi(|\zeta_1-\zeta_2|)\big]^{\frac{2-p}{2}} \left[\Psi_{|\zeta_2|}(|\zeta_1-\zeta_2|)\right]^{\frac{p}{2}}.
\end{align}
For every $\epsilon \in (0,1)$, let us apply Young's inequality on the right-hand side of~\eqref{est-phi-PQ}, it follows that
\begin{align}\notag
\Psi(|\zeta_1-\zeta_2|) &\le \frac{\epsilon}{2} \big[\Psi(|\zeta_1|) + \Psi(|\zeta_1-\zeta_2|)\big] + C \epsilon^{1-\frac{2}{p}} \Psi_{|\zeta_2|}(|\zeta_1-\zeta_2|) \\
&\le \frac{1}{2} \Psi(|\zeta_1-\zeta_2|) + \frac{\epsilon}{2} \Psi(|\zeta_1|) + C \epsilon^{1-\frac{2}{p}} \Psi_{|\zeta_2|}(|\zeta_1-\zeta_2|),\notag
\end{align}
which allows us to conclude~\eqref{vphi-2} by combining with~\eqref{Phi-PQ}.
\end{proof}

\begin{lemma}\label{lem:comp}
Suppose that $u \in \mathsf{g} + W^{1,p}_{0,\omega}(\Omega)$ is a weak solution to~\eqref{Eq-0} under assumptions in Lemma~\ref{lem:global}. Let $x_0 \in \overline{\Omega}$ and $B = B(x_0,R)$, we denote 
$$\lambda B = B(x_0,\lambda R) \ \mbox{ and } \ \lambda\Omega_{B} = \lambda B \cap \Omega \ \mbox{ for } \lambda>0.$$
There exists a function $v \in W^{1,p}_{\omega}(\Omega_B)$ such that for every $\epsilon \in (0,1)$, there holds
\begin{align}\label{est-lem:comp}
\fint_{\Omega_B} |\mathbb{P}(x)\nabla u & - \mathbb{P}(x)\nabla v|^p  dx  \le \epsilon \fint_{\Omega_B} |\mathbb{P}(x)\nabla u|^{p}dx \notag \\
& \qquad   + C_{\epsilon} \left(\fint_{\Omega_B} |\mathbb{P}(x)\mathbf{F}|^{p}dx + \fint_{\Omega_B} |\mathbb{P}(x)\nabla \mathsf{g}|^{p}dx\right).
\end{align}
Moreover, there exists a constant $\kappa>0$ such that if $|\log \mathbb{P}|_{\mathrm{BMO}} \le \kappa$ and $\Omega$ is $(\kappa,r_0)$-Lipschitz for some $r_0>0$, then the following inequality
\begin{align}\label{est-lem:Re}
\left(\fint_{\frac{1}{16} \Omega_B} |\mathbb{P}(x)\nabla v|^{p\gamma}dx  \right)^{\frac{1}{\gamma}} \le C(\gamma) \left[ \fint_{\frac{1}{2} \Omega_B} |\mathbb{P}(x)\nabla v|^{p} dx + \left(\fint_{\frac{1}{2} \Omega_B} |\mathbb{P}(x)\nabla \mathsf{g}|^{p\gamma}dx \right)^{\frac{1}{\gamma}}\right]
\end{align}
holds for every $\gamma \ge 1$. 
\end{lemma}
\begin{proof}
Let $v \in u - \mathsf{g} + W^{1,p}_{0,\omega}(\Omega_B)$ be the weak solution to the following problem  
\begin{align}\label{Eq-v}
\begin{cases}
-\mathrm{div}\left(\mathcal{L}_p(x,\nabla v)\right) & = \ 0  \quad \mbox{in} \ \Omega_B,\\
\hspace{1cm} v &  = \ u - \mathsf{g}  \quad \mbox{on} \ \partial \Omega_B.\end{cases}
\end{align}
Therefore, $v$ solves the following variational formula
\begin{align}\label{var-v}
& \int_{\Omega_B} |\mathbb{P}(x)\nabla v|^{p-2}\mathbb{P}(x)\nabla v \cdot  \mathbb{P}(x)\nabla\varphi dx  =  0
\end{align}
for all $\varphi \in W^{1,p}_{0,\omega}(\Omega_B)$. Testing~\eqref{var-P0} and~\eqref{var-v} by $\varphi = u - v - \mathsf{g}$, we obtain that
\begin{align}
\fint_{\Omega_B} & \left(|\mathbb{P}(x)\nabla u|^{p-2}\mathbb{P}(x)\nabla u - |\mathbb{P}(x)\nabla v|^{p-2}\mathbb{P}(x)\nabla v\right) \cdot \left(\mathbb{P}(x)\nabla u - \mathbb{P}(x) \nabla v\right)  dx \le J, \label{pi-u-v-1}
\end{align}
where $J$ is given by
\begin{align*}
J & := \fint_{\Omega_B} |\mathbb{P}(x)\nabla u - \mathbb{P}(x) \nabla v| |\mathbb{P}(x)\mathbf{F}|^{p-1}  dx + \fint_{\Omega_B} |\mathbb{P}(x)\nabla \mathsf{g}| |\mathbb{P}(x)\nabla u|^{p-1}  dx \notag \\
& \qquad  + \fint_{\Omega_B} |\mathbb{P}(x)\nabla \mathsf{g}| |\mathbb{P}(x)\nabla v|^{p-1}  dx + \fint_{\Omega_B} |\mathbb{P}(x)\nabla \mathsf{g}| |\mathbb{P}(x)\mathbf{F}|^{p-1}  dx.
\end{align*}
Using the notation in~\eqref{def-Phi}, it is well-known that
\begin{align*}
|\mathbb{V}_{p}(\zeta_1)-\mathbb{V}_{p}(\zeta_2)|^2 \simeq \left(|\zeta_1|^{p-2}\zeta_1 - |\zeta_2|^{p-2}\zeta_2\right) \cdot (\zeta_1-\zeta_2), \ \mbox{ for all } \zeta_1, \zeta_2 \in \mathbb{R}^n.
\end{align*}
Hence,~\eqref{pi-u-v-1} implies to
\begin{align}\label{pi-u-v-2}
\fint_{\Omega_B} |\mathbb{V}_{p}(\mathbb{P}(x) \nabla u) & - \mathbb{V}_{p}(\mathbb{P}(x) \nabla v)|^2 dx \le CJ.
\end{align}
For $p \ge 2$, thanks to~\eqref{vphi-1} in Lemma~\ref{lem-Phi} and~\eqref{pi-u-v-2}, one gets that
\begin{align}\label{pi-u-v-3}
\fint_{\Omega_B} |\mathbb{P}(x)\nabla u & - \mathbb{P}(x)\nabla v|^p dx  \le C J.
\end{align}
Applying Young's inequality for all terms of $J$, it implies to~\eqref{est-lem:comp} from~\eqref{pi-u-v-3}. For $1<p<2$, we will apply~\eqref{vphi-2} in Lemma~\ref{lem-Phi}, it follows from~\eqref{pi-u-v-2} that
\begin{align}
\fint_{\Omega_B} |\mathbb{P}(x)\nabla u - \mathbb{P}(x)\nabla v|^p dx & \le \epsilon \fint_{\Omega_B} |\mathbb{P}(x)\nabla u|^p dx \notag \\
& \qquad + C {\epsilon}^{1 -\frac{2}{p}} \fint_{\Omega_B} |\mathbb{V}_{p}(\mathbb{P}(x) \nabla u) - \mathbb{V}_{p}(\mathbb{P}(x) \nabla v)|^2 dx. \notag
\end{align}
It yields to~\eqref{est-lem:comp} by combining with~\eqref{pi-u-v-2} and  using Young's inequality for the last integral term.

The reverse H\"older's inequality~\eqref{est-lem:Re} is a consequence of the main results in~\cite{BDGN2022} and~\cite{BBDL2023} for the homogeneous problem~\eqref{Eq-v}. More precisely, if $B \subset \Omega$ then by~\cite[Theorem 2]{BDGN2022}, inequality~\eqref{est-lem:Re} holds provided 
$$|\log \mathbb{P}|_{\mathrm{BMO}} \le \kappa, \quad \mbox{ for } \kappa \mbox{ small enough.}$$ 
Otherwise, if $B \cap \partial \Omega \neq \emptyset$ then~\eqref{est-lem:Re} is deduced from inequality (3.123) in~\cite{BBDL2023}. In this boundary case, an additional assumption that $\Omega$ is $(\kappa,r_0)$-Lipschitz for some $r_0>0$, are made. The proof is complete.
\end{proof}

A large-scaling property of level-set inequality will be implemented based on the covering lemma. It nowadays becomes standard in the argument of several approaches in the literature. For convenience, we restate here a modified version of Calder\'on-Zygmund covering lemma as below, the interested reader may consult~\cite{CC1995, Caffa1998, Vitali}. 

\begin{lemma}\label{lem:Vitali}
Let $\Omega$ be a $(\kappa, r_0)$-Lipschitz domain with $\kappa, r_0>0$ (Assumption $A_2$). Suppose that $\mu \in \mathrm{A}_{\infty}$ and two measurable subsets $D\subset E$ of $\Omega$ satisfy:
\begin{enumerate}
\item[i)] $\mu(D) < \varepsilon \mu(B(0,R_0))$, for some $\varepsilon \in (0,1)$ and $0<R_0\le r_0$;
\item[ii)] if $B(x,\rho)\cap \Omega \not\subset E$ then $\mu(B(x,\rho)\cap D) < \varepsilon \mu(B(x,\rho))$, for every $x\in \Omega$ and $0<\rho \le R_0$.
\end{enumerate} 	
Then, there exists a constant $C>0$ only depending on $n$ such that $\mu(D) \leq C\varepsilon \mu(E)$.
\end{lemma}

\begin{lemma}\label{lem-Sig}
Suppose that the non-decreasing function $\Sigma: [0,\infty) \to [0,\infty)$  satisfies the following doubling property
\begin{align}\label{cond-I}
c_1 \Sigma(t) \le \Sigma(2t) \le c_2 \Sigma(t), \mbox{ for all } t \ge 0
\end{align}
for two constants $c_1, c_2>1$. Then, there holds
\begin{align}\label{ineq-Sig-1}
\Sigma(\sigma_1 + \sigma_2) \le c_2 \big[\Sigma(\sigma_1) + \Sigma(\sigma_2)\big], \mbox{ for all } \sigma_1, \sigma_2 \ge 0.
\end{align}
Moreover, for every $\epsilon \in (0,1/2)$ and $t \ge 0$, there holds
\begin{align}\label{ineq-Sig-2}
\Sigma(\epsilon t) \le c_1\epsilon^{\log_2 c_1} \Sigma(t).
\end{align}
\end{lemma}
\begin{proof}
The proof of~\eqref{ineq-Sig-1} is very simple. Indeed, combining the second inequality in~\eqref{cond-I} and the fact that $\Sigma$ is non-deceasing, one has
\begin{align}\notag
\Sigma(\sigma_1 + \sigma_2) \le \Sigma(2\max\{\sigma_1,\sigma_2\}) \le c_2 \Sigma(2\max\{\sigma_1,\sigma_2\}) \le c_2 \big[\Sigma(\sigma_1) + \Sigma(\sigma_2)\big], 
\end{align}
for all $\sigma_1, \sigma_2 \ge 0$. Let us now prove~\eqref{ineq-Sig-2}. For every $\epsilon \in (0,1/2)$, one can find $k = k(\epsilon) \in \mathbb{Z}^+$ satisfying
\begin{align*}
2^{-k-1} < \epsilon \le 2^{-k} \Leftrightarrow \log_2\epsilon \le -k < 1+\log_2\epsilon.
\end{align*}
Applying the first inequality in~\eqref{cond-I}, there holds
\begin{align*}
\Sigma(\epsilon t) \le \Sigma(2^{-k} t) \le c_1^{-k} \Sigma(t) \le c_1^{1+\log_2\epsilon} \Sigma(t),
\end{align*}
which leads to~\eqref{ineq-Sig-2}. 
\end{proof}

\begin{lemma}\label{lem:Mor}
Let $y \in \Omega$, $0<\varrho<D_0$. Then the following estimate holds
\begin{align}\label{gmu}
\left[\mathbf{M}\chi_{B(y,\varrho)}\right](x) \le 2^{-(j-1)n}, \quad \mbox{ for all } x \in \Omega_j^{\varrho}(y),
\end{align}
where $\Omega_j$ is defined by
\begin{align}
\Omega_j^{\varrho}(y) = \left\{x\in \mathbb{R}^n: \ 2^{j}\varrho \le |x-y| < 2^{j+1}\varrho\right\}, \quad j \in \mathbb{Z}^+.\label{Omega-j}
\end{align}
\end{lemma}
\begin{proof}
Let $j \in \mathbb{Z}^+$ and $x \in \Omega_j^{\varrho}(y)$. One has 
\begin{align}\notag
\left[\mathbf{M}\chi_{B(y,\varrho)}\right](x) & = \sup_{r>0} \fint_{B(x,r)} \chi_{B(y,\varrho)}(z) dz =  \sup_{r>0} \mathcal{R}((y,\varrho);(x,r)),
\end{align}
where the ratio $\mathcal{R}((y,\varrho);(x,r))$ is defined by
\begin{align*}
\mathcal{R}((y,\varrho);(x,r)) := \frac{|B(y,\varrho) \cap B(x,r)|}{|B(x,r)|}.
\end{align*}
Since $2^{j}\varrho \le |x - y| < 2^{j+1}\varrho$, there holds 
$$|z - x| \ge |x-y| - |y - z| > 2^j\varrho-\varrho \ge 2^{j-1}\varrho,$$
for every $z \in B(y,\varrho)$. It implies to 
$$|B(y,\varrho)\cap B(x,r)| = 0 \ \mbox{ for all } \ r \le 2^{j-1}\varrho.$$ 
Similarly, one can check that $B(x,r) \supset B(y,\varrho)$ for every $r \ge 2^{j+2}\varrho$, thus
\begin{align*}
\sup_{r \ge 2^{j+2}\varrho} \mathcal{R}((y,\varrho);(x,r)) = \sup_{r \ge 2^{j+2}\varrho} {|B(x,r)|^{-1}} {|B(y,\varrho)|} = 2^{-(j+2)n}. 
\end{align*}
On the other hand, one has
\begin{align}\notag
\sup_{2^{j-1}\varrho < r < 2^{j+2}\varrho} \mathcal{R}((y,\varrho);(x,r)) &\le \sup_{2^{j-1}\varrho < r < 2^{j+2}\varrho} {|B(x,r)|^{-1}} {|B(y,\varrho)|}  = 2^{-(j-1)n}.\notag
\end{align}
Taking into account all above estimates, one may conclude~\eqref{gmu}.
\end{proof}

\section{Proofs of main theorems}
\label{sec:proof}

We are now ready to prove our main results. It is worth noticing that regarding Theorem~\ref{theo:general}, we shall split the statement into some small theorems associated with each subtle function space introduced in Section~\ref{sec:pre}.

\begin{proof}[Proof of Theorem~\ref{theo:A}]
Let us first introduce two subsets
\begin{align}\notag 
\mathbb{D} & := \left\{\mathbf{M}_{\alpha}(|\mathbb{P}(x)\nabla u|^{p}) > \varepsilon^{-\theta} \lambda; \, \mathbf{M}_{\alpha}\big(|\mathbb{P}(x) \mathbf{G}|^{p}\big) \le \varepsilon^{\gamma} \lambda\right\},\\
\mathbb{E} & := \{\mathbf{M}_{\alpha}(|\mathbb{P}(x)\nabla u|^{p}) > \lambda\}.\notag
\end{align}
The following decomposition
\begin{align*}
\left\{\mathbf{M}_{\alpha}(|\mathbb{P}(x)\nabla u|^{p}) > \varepsilon^{-\theta} \lambda\right\} = \mathbb{D} \cup \left\{\mathbf{M}_{\alpha}(|\mathbb{P}(x)\nabla u|^{p}) > \varepsilon^{-\theta} \lambda; \, \mathbf{M}_{\alpha}\big(|\mathbb{P}(x) \mathbf{G}|^{p}\big) > \varepsilon^{\gamma} \lambda\right\}
\end{align*}
allows us to arrive
\begin{align}\label{DVF}
\mu\left(\left\{\mathbf{M}_{\alpha}(|\mathbb{P}(x)\nabla u|^{p}) > \varepsilon^{-\theta} \lambda\right\}\right) \le \mu(\mathbb{D}) + \mu\left(\left\{\mathbf{M}_{\alpha}\big(|\mathbb{P}(x) \mathbf{G}|^{p}\big) > \varepsilon^{\gamma} \lambda\right\}\right).
\end{align}
If the following inequality holds
\begin{align}\label{goal-1}
\mu(\mathbb{D}) \le C \varepsilon \mu(\mathbb{E}),
\end{align}
then~\eqref{DVF} implies to~\eqref{dist-ineq}. For this reason, it sufficient to prove~\eqref{goal-1}. Thanks to Lemma~\ref{lem:Vitali}, we will show two statements:
\begin{enumerate}
\item[i)] $\mu(\mathbb{D}) < \varepsilon \mu(B(0,R_0))$ for $0<R_0\le r_0$;
\item[ii)] if $B(x,\rho)\cap \Omega \not\subset \mathbb{E}$ then $\mu(B(x,\rho)\cap \mathbb{D}) < \varepsilon \mu(B(x,\rho))$, for every $x\in \Omega$ and $0<\rho \le R_0$.
\end{enumerate}

The first statement $i)$ is valid if $\mathbb{D}$ is empty. Otherwise, one can find $\zeta_1 \in \Omega$ such that $\mathbf{M}_{\alpha}\big(|\mathbb{P} \mathbf{G}|^{p}\big)(\zeta_1) \le \varepsilon^{\gamma} \lambda$, which leads to
\begin{align}\label{cond-x1}
\fint_{B(\zeta_1,\varrho)} |\mathbb{P}(x) \mathbf{G}|^{p} dx \le \varepsilon^{\gamma} \varrho^{-\alpha} \lambda, \mbox{ for all } \varrho>0.
\end{align}
Assume that $[\mu]_{\mathrm{A}_\infty} = (c_1,c_2,\nu_1,\nu_2)$. It is possible to find a ball $B(\zeta_1, R_1)$ such that
$$\Omega \cup B(0,R_0) \subset B(\zeta_1, R_1).$$
We remark that the ratio $R_1/R_0$ depends on $D_0/R_0$. Thanks to~\eqref{Muck-w}, since $\mathbb{D} \subset \Omega \subset B(\zeta_1, R_1)$,  one has
\begin{align}
\mu(\mathbb{D}) &\le c_2\left(\frac{|\mathbb{D}|}{|B(\zeta_1, R_1)|}\right)^{\nu_2}\mu(B(\zeta_1, R_1)) \notag \\
&\le c_2\left(\frac{|\mathbb{D}|}{|B(\zeta_1, R_1)|}\right)^{\nu_2}c_1^{-1}\left(\frac{|B(\zeta_1, R_1)|}{|B(0,R_0)|}\right)^{\nu_1}\mu(B(0,R_0)) \notag \\
& \le C\left(\frac{|\mathbb{D}|}{|B(\zeta_1, R_1)|}\right)^{\nu_2} \mu(B(0,R_0)).\label{est-i-1a}
\end{align}
Thanks to Lemma~\ref{lem:bound-M-beta} and inequality~\eqref{est-global} in Lemma~\ref{lem:global}, there holds
\begin{align}
|\mathbb{D}| & \le \left|\left\{\mathbf{M}_{\alpha}(|\mathbb{P}(x)\nabla u|^{p}) > \varepsilon^{-\theta} \lambda\right\}\right| \notag \\
& \le C\left(\varepsilon^{\theta} \lambda^{-1}\int_{\Omega} |\mathbb{P}(x)\nabla u|^{p}  dx\right)^{\frac{n}{n-\alpha}} \notag \\
& \le C\left(\varepsilon^{\theta} \lambda^{-1}\int_{\Omega} |\mathbb{P}(x) \mathbf{G}|^{p}  dx\right)^{\frac{n}{n-\alpha}} \notag \\
& \le C\left[\varepsilon^{\theta} \lambda^{-1} |B(\zeta_1, R_1)| \fint_{B(\zeta_1, R_1)} |\mathbb{P}(x) \mathbf{G}|^{p}  dx\right]^{\frac{n}{n-\alpha}}. \label{est-i-1}
\end{align}
Substituting~\eqref{cond-x1} into~\eqref{est-i-1}, it yields that
\begin{align}
|\mathbb{D}| & \le C\left[\varepsilon^{\theta} \lambda^{-1} |B(\zeta_1, R_1)|^{1-\frac{\alpha}{n}}  \varepsilon^{\gamma} \lambda \right]^{\frac{n}{n-\alpha}} \le C \varepsilon^{(\theta+\gamma)\frac{n}{n-\alpha}} |B(\zeta_1, R_1)|. \label{est-i-2}
\end{align}
Combining~\eqref{est-i-1a} and~\eqref{est-i-2}, it follows that
\begin{align}
\mu(\mathbb{D}) \le C \varepsilon^{(\theta+\gamma)\frac{n\nu_2}{n-\alpha}} \mu(B(0,R_0)) < \varepsilon \mu(B(0,R_0)),\notag
\end{align}
for every $\varepsilon$ small enough, which satisfies  $C \varepsilon^{(\theta+\gamma)\frac{n\nu_2}{n-\alpha}} < \varepsilon$ if $(\theta+\gamma)\frac{n\nu_2}{n-\alpha}>1$.\\

Let us now prove $ii)$. Let $\zeta\in \Omega$ and $\rho \in (0,R_0]$ such that $B(\zeta,\rho)\cap \Omega \not\subset \mathbb{E}$, we will show that
\begin{align}\label{goal-2}
\mu(\mathbb{D}\cap B(\zeta,\rho)) < \varepsilon \mu(B(\zeta,\rho)).
\end{align}
From now on, we will denote $B = B(\zeta,\rho)$ and $rB = B(\zeta, r\rho)$ for simplicity. By assuming $\mathbb{D}\cap B \neq \emptyset$, one can find  $\zeta_2, \zeta_3 \in B \cap \Omega$ such that $\mathbf{M}_{\alpha}(|\mathbb{P} \nabla u|^{p})(\zeta_2) \le \lambda$ and $\mathbf{M}_{\alpha}\big(|\mathbb{P} \mathbf{G}|^{p}\big)(\zeta_3) \le \varepsilon^{\gamma}\lambda$. It implies to
\begin{align}\label{cond-x23}
\fint_{B(\zeta_2,\varrho)} |\mathbb{P}(x) \nabla u|^{p} dx \le \varrho^{-\alpha}\lambda, \mbox{ and }
 \fint_{B(\zeta_3,\varrho)} |\mathbb{P}(x) \mathbf{G}|^{p} dx \le \varepsilon^{\gamma} \varrho^{-\alpha}  \lambda, \quad \forall \varrho>0.
\end{align}
For every $y \in B$, one can check that $B(y,\varrho) \subset B(\zeta_2,3\varrho)$ for all $\varrho \ge \rho$. Indeed, for each $z \in B(y,\varrho)$, it follows that
$$|z - \zeta_2| \le |z-y| + |y -\zeta| + |\zeta - \zeta_2| < \varrho + \rho + \rho \le 3 \varrho.$$
For this reason, by~\eqref{cond-x23}, one has\begin{align*}
\mathbf{M}_{\alpha}(|\mathbb{P}\nabla u|^{p})(y) & = \max\left\{\sup_{0<\varrho<\rho} \varrho^{\alpha} \fint_{B(y,\varrho)} |\mathbb{P}(x)\nabla u|^{p} dx; \, \sup_{\varrho \ge \rho} \varrho^{\alpha} \fint_{B(y,\varrho)} |\mathbb{P}(x)\nabla u|^{p} dx\right\} \\
& \le \max\left\{\mathbf{M}_{\alpha}^{\rho}(\chi_{2B}|\mathbb{P}\nabla u|^{p})(y); \,  3^n \sup_{\varrho \ge \rho} \varrho^{\alpha} \fint_{B(\zeta_2,3\varrho)} |\mathbb{P}(x)\nabla u|^{p} dx\right\} \\
& \le \max\left\{\mathbf{M}_{\alpha}^{\rho}(\chi_{2B}|\mathbb{P}\nabla u|^{p})(y); \,  3^n \lambda\right\}.
\end{align*}
Therefore, for $\varepsilon^{-\theta}>3^n$, it holds
\begin{align}\label{est-ii-1}
\mathbb{D}\cap B \subset \left\{\mathbf{M}_{\alpha}^{\rho}(\chi_{2B}|\mathbb{P}(x)\nabla u|^{p}) > \varepsilon^{-\theta}\lambda\right\} \cap B.
\end{align}
If $\zeta$ is far from the boundary of $\Omega$, we may assume that $4B \subset \Omega$. Otherwise, we assume that $4B \cap \partial \Omega \neq \emptyset$. In both cases, we may cover $2B$ by a new ball $\tilde{B}(\tilde{\zeta},\tilde{\rho})$ for $\tilde{\zeta} \in \overline{\Omega}$. Indeed, if $4B \subset \Omega$, we take $\tilde{\zeta} = \zeta$ and $\tilde{\rho} = 4\rho$. If $4B \cap \partial \Omega \neq \emptyset$, we take $\tilde{\rho} = 6\rho$ and $\tilde{\zeta} \in 4B \cap \partial \Omega$ such that 
$|\tilde{\zeta}-\zeta| = d(\zeta, \partial \Omega)$. We now denote
\begin{align*}
r\Omega_{\tilde{B}} := \Omega \cap B(\tilde{\zeta},r\tilde{\rho}), \quad \mbox{ for } r>0.
\end{align*}
Then, we rewrite~\eqref{est-ii-1} as follows
\begin{align}\label{est-ii-1c}
\mathbb{D}\cap B \subset \left\{\mathbf{M}_{\alpha}^{\rho}(\chi_{\Omega_{\tilde{B}}}|\mathbb{P}(x)\nabla u|^{p}) > \varepsilon^{-\theta}\lambda\right\} \cap B.
\end{align}
Thanks to Lemma~\ref{lem:comp}, there exist $v \in W^{1,p}_{\omega}(8\Omega_{\tilde{B}})$ and $\kappa>0$ such that if $|\log \mathbb{P}|_{\mathrm{BMO}} \le \kappa$ and $\Omega$ is $(\kappa,r_0)$-Lipschitz  then
\begin{align}\label{est-ii-2}
\left(\fint_{\Omega_{\tilde{B}}} |\mathbb{P}(x)\nabla v|^{p\gamma}dx\right)^{\frac{1}{\gamma}} \le C \left[\fint_{8\Omega_{\tilde{B}}} |\mathbb{P}(x)\nabla v|^{p} dx + \left(\fint_{\Omega_{8\tilde{B}}} |\mathbb{P}(x)\nabla \mathsf{g}|^{p\gamma}dx\right)^{\frac{1}{\gamma}}\right],
\end{align}
for all $\gamma \ge 1$ and
\begin{align}\label{est-ii-3}
\fint_{8\Omega_{\tilde{B}}} |\mathbb{P}(x)\nabla u - \mathbb{P}(x)\nabla v|^p dx & \le \epsilon \fint_{8\Omega_{\tilde{B}}} |\mathbb{P}(x)\nabla u|^{p}dx   + C_{\epsilon} \fint_{8\Omega_{\tilde{B}}} |\mathbb{P}(x)\mathbf{G}|^{p}dx,
\end{align}
for every $\epsilon \in (0,1)$. Since $\tilde{\rho} \le 6\rho$, one has
$$8\tilde{B} \subset 48 B(\tilde{\zeta},\rho) \subset 52 B \subset 53 B(\zeta_2,\rho) \cap 53 B(\zeta_3,\rho),$$ 
which by~\eqref{cond-x23} ensures that 
\begin{align}\label{est-piu}
\fint_{8\Omega_{\tilde{B}}} |\mathbb{P}(x)\nabla u|^{p}dx & \le C \fint_{53 B(\zeta_2,\rho)} |\mathbb{P}(x)\nabla u|^{p}dx \le C \rho^{-\alpha} \lambda,
\end{align}
and
\begin{align}\label{est-piF}
\fint_{8\Omega_{\tilde{B}}} |\mathbb{P}(x) \mathbf{G}|^{p}dx & \le C \fint_{53 B(\zeta_3,\rho)} |\mathbb{P}(x)\mathbf{G}|^{p}dx \le C \rho^{-\alpha} \varepsilon^{\gamma} \lambda.
\end{align}
Substituting~\eqref{est-piu} and~\eqref{est-piF} into~\eqref{est-ii-3}, one gets that
\begin{align}
 \fint_{8\Omega_{\tilde{B}}} |\mathbb{P}(x)\nabla u - \mathbb{P}(x)\nabla v|^p dx  & \le C \rho^{-\alpha} \left(\epsilon + C_{\epsilon} \varepsilon^{\gamma}\right) \lambda.\notag 
\end{align}
For simplicity of computation, we may fix $\epsilon$ and $\gamma$ such that 
$\epsilon = C_{\epsilon} \varepsilon^{\gamma} = \varepsilon^{\delta}$,
for a new positive exponent $\delta$ determined later. It follows that
\begin{align}
 \fint_{8\Omega_{\tilde{B}}} |\mathbb{P}(x)\nabla u - \mathbb{P}(x)\nabla v|^p dx  & \le C \rho^{-\alpha} \varepsilon^{\delta} \lambda.\label{est-ii-3b} 
\end{align}
Using a simple inequality on the right-hand side of~\eqref{est-ii-2} and taking~\eqref{est-piu},~\eqref{est-ii-3b} into account, one has
\begin{align}
\fint_{\Omega_{\tilde{B}}} |\mathbb{P}(x)\nabla v|^{p\gamma}dx  &\le C \left(\fint_{8\Omega_{\tilde{B}}} |\mathbb{P}(x)\nabla u|^{p} dx\right)^{\gamma} \notag \\
& \qquad + C \left(\fint_{8\Omega_{\tilde{B}}} |\mathbb{P}(x)\nabla u - \mathbb{P}(x)\nabla v|^{p} + |\mathbb{P}(x)\nabla \mathsf{g}|^{p} dx \right)^{\gamma}\notag \\
& \le C \rho^{-\alpha\gamma} \left(1 + \varepsilon^{\delta\gamma}\right) \lambda^{\gamma} \notag \\
& \le C \rho^{-\alpha\gamma} \lambda^{\gamma}. \label{est-ii-2b}
\end{align}
Using a fundamental inequality, from~\eqref{est-ii-1c}, one gets that
\begin{align}
|\mathbb{D}\cap B| &\le \left|\left\{\mathbf{M}_{\alpha}^{\rho}(\chi_{\Omega_{\tilde{B}}}|\mathbb{P}(x)\nabla v|^{p}) > 2^{-p} \varepsilon^{-\theta}\lambda\right\}\right|  \notag \\
& \qquad \qquad + \left|\left\{\mathbf{M}_{\alpha}^{\rho}(\chi_{\Omega_{\tilde{B}}}|\mathbb{P}(x)\nabla u - \mathbb{P}(x)\nabla v|^{p}) > 2^{-p}\varepsilon^{-\theta}\lambda\right\}\right|.\notag
\end{align}
Thanks to Lemma~\ref{lem:bound-M-beta}, by~\eqref{est-ii-3b} and~\eqref{est-ii-2b} it deduces that
\begin{align}
|\mathbb{D}\cap B| & \le C \left(2^{p\gamma} \varepsilon^{\theta\gamma}\lambda^{-\gamma} \rho^n \fint_{\Omega_{\tilde{B}}}|\mathbb{P}(x)\nabla v|^{p\gamma}dx\right)^{\frac{n}{n-\alpha\gamma}} \notag \\
& \qquad + C \left(2^{p} \varepsilon^{\theta}\lambda^{-1} \rho^n \fint_{\Omega_{\tilde{B}}}|\mathbb{P}(x)\nabla u - \mathbb{P}(x)\nabla v|^{p}dx\right)^{\frac{n}{n-\alpha}}\notag \\
 & \le C \left(\varepsilon^{\theta\gamma} \rho^{n-\alpha\gamma}\right)^{\frac{n}{n-\alpha\gamma}}  + C \left(\varepsilon^{\theta+\delta} \rho^{n-\alpha}\right)^{\frac{n}{n-\alpha}}\notag \\
& \le C \left(\varepsilon^{\frac{\theta n\gamma}{n-\alpha\gamma}} + \varepsilon^{\frac{n(\theta+\delta)}{n-\alpha}}\right) |B|.\notag
\end{align}
Applying~\eqref{Muck-w} again, it yields that
\begin{align}
\mu(\mathbb{D}\cap B) & \le C \left(\varepsilon^{\frac{\theta n\gamma\nu_2}{n-\alpha\gamma}} + \varepsilon^{\frac{n(\theta+\delta)\nu_2}{n-\alpha}}\right) \mu(B).\label{est-ii-5}
\end{align}
One can see that if the exponents of $\varepsilon$ are larger than 1 then~\eqref{est-ii-5} implies to~\eqref{goal-2} for every $\varepsilon$ small enough. To do this, we simply choose suitable values of $\gamma>1$ and $\delta>0$ large enough such that
\begin{align}\notag 
\min\left\{\frac{\theta n\gamma\nu_2}{n-\alpha\gamma}; \, \frac{n(\theta+\delta)\nu_2}{n-\alpha}\right\} > 1.
\end{align}
With these choices, there exists $\varepsilon$ small enough such that
\begin{align*}
C \left(\varepsilon^{\frac{\theta n\gamma\nu_2}{n-\alpha\gamma}} + \varepsilon^{\frac{n(\theta+\delta)\nu_2}{n-\alpha}}\right) < \varepsilon,
\end{align*}
which completes the proof of~\eqref{goal-2} from~\eqref{est-ii-5}. 
\end{proof}

Our first application of the general level-set argument concerns the weighted Lorentz spaces.

\begin{theorem}\label{theo:B}
Let $\mathbb{P}$ be a matrix weight satisfying~\eqref{cond-PI}; $\mathbf{F} \in L^{p}_{\omega}(\Omega)$ and $\mathsf{g} \in W^{1,p}_{\omega}(\overline{\Omega})$ for $p \in (1,\infty)$. Assume further that $u \in \mathsf{g} + W^{1,p}_{0,\omega}(\Omega)$ is a weak solution to~\eqref{Eq-0}, a Muckenhoupt weight $\mu \in \mathrm{A}_{\infty}$, $\alpha \in [0,n)$. Let $0<\mathfrak{q}<\infty$ and $0<\mathfrak{s}\le\infty$, then there exists $\kappa = \kappa(\mathtt{data},\mathfrak{q},\mathfrak{s})>0$ such that 
\begin{align}\label{est:B}
\|\mathbf{M}_{\alpha}(|\mathbb{P}\nabla u|^{p})\|_{L^{\mathfrak{q},\mathfrak{s}}_{\mu}(\Omega)}  \le C\left\|\mathbf{M}_{\alpha}(|\mathbb{P}\mathbf{G}|^{p})\right\|_{L^{\mathfrak{q},\mathfrak{s}}_{\mu}(\Omega)},
\end{align}
 if $\log \mathbb{P}$ satisfies $(\kappa,r_0)$-small-$\log$-$\mathrm{BMO}$ condition and $\Omega$ is $(\kappa,r_0)$-Lipschitz domain for some $r_0>0$. Here, the $C$ is a positive constant depending on $\mathfrak{q},\mathfrak{s},\mathtt{data}$.
\end{theorem}

\begin{proof}[Proof of Theorem~\ref{theo:B}]
Let $0<\mathfrak{q}<\infty$ and $0<\mathfrak{s}<\infty$. Applying Theorem~\ref{theo:A} with $\theta = \frac{1}{2\mathfrak{q}}$, there exist positive constants 
$$\gamma = \gamma(\alpha,\mathfrak{q},\mathfrak{s})>0 \ \mbox{ and } \ \kappa = \kappa(\mathfrak{q},\mathfrak{s},\varepsilon,\alpha)>0$$ 
such that if $|\log \mathbb{P}|_{\mathrm{BMO}} \le \kappa$ and $\Omega$ is $(\kappa,r_0)$-Lipschitz for some $r_0>0$, then the following level-set inequality
\begin{align}\notag
 & d^{\mu}_{\alpha}\big(|\mathbb{P}(x)\nabla u|^{p};\varepsilon^{-\frac{1}{2\mathfrak{q}}}\lambda\big) \le C \varepsilon d^{\mu}_{\alpha}(|\mathbb{P}(x)\nabla u|^{p};\lambda) + d^{\mu}_{\alpha}(|\mathbb{P}(x)\mathbf{G}|^{p};\varepsilon^{\gamma}\lambda)
\end{align}
holds for any $\lambda> 0$ and $\varepsilon$ small enough. Replacing $\varepsilon^{-\frac{1}{2\mathfrak{q}}}\lambda$ by $\lambda$, one rewrite this inequality as 
\begin{align}\notag
 & d^{\mu}_{\alpha}\big(|\mathbb{P}(x)\nabla u|^{p};\lambda\big) \le C \varepsilon d^{\mu}_{\alpha}\big(|\mathbb{P}(x)\nabla u|^{p};\varepsilon^{\frac{1}{2\mathfrak{q}}}\lambda\big) + d^{\mu}_{\alpha}\big(|\mathbb{P}(x)\mathbf{G}|^{p};\varepsilon^{\gamma+\frac{1}{2\mathfrak{q}}}\lambda\big).
\end{align}
Multiplying by $\lambda^{\mathfrak{q}}$ both sides of the above inequality and then taking the integral over $\mathbb{R}^+$ with respect to $\lambda$, one obtains that
\begin{align}\notag
 & \int_{\mathbb{R}^+} \mathfrak{q}\left[\lambda^{\mathfrak{q}} d^{\mu}_{\alpha}\big(|\mathbb{P}(x)\nabla u|^{p};\lambda\big) \right]^{\frac{\mathfrak{s}}{\mathfrak{q}}} \frac{d\lambda}{\lambda}  \le C \varepsilon^{\frac{\mathfrak{s}}{\mathfrak{q}}}  \int_{\mathbb{R}^+}\mathfrak{q} \left[\lambda^{\mathfrak{q}} d^{\mu}_{\alpha}\big(|\mathbb{P}(x)\nabla u|^{p};\varepsilon^{\frac{1}{2\mathfrak{q}}}\lambda\big)\right]^{\frac{\mathfrak{s}}{\mathfrak{q}}} \frac{d\lambda}{\lambda} \\
 & \hspace{2cm} + C \int_{\mathbb{R}^+} \mathfrak{q}\left[\lambda^{\mathfrak{q}} d^{\mu}_{\alpha}\big(|\mathbb{P}(x)\mathbf{G}|^{p};\varepsilon^{\gamma+\frac{1}{2\mathfrak{q}}}\lambda\big)\right]^{\frac{\mathfrak{s}}{\mathfrak{q}}} \frac{d\lambda}{\lambda}.\label{LR-1}
\end{align}
Changing of variables for two terms on the right-hand side of~\eqref{LR-1} and using the quasi-norm in~\eqref{w-Lor-norm}, it deduces to 
\begin{align}\notag
 \|\mathbf{M}_{\alpha}(|\mathbb{P}(x)\nabla u|^{p})\|_{L^{\mathfrak{q},\mathfrak{s}}_{\mu}(\Omega)}^{\mathfrak{s}} & \le C \varepsilon^{\frac{\mathfrak{s}}{2\mathfrak{q}}}  \|\mathbf{M}_{\alpha}(|\mathbb{P}(x)\nabla u|^{p})\|_{L^{\mathfrak{q},\mathfrak{s}}_{\mu}(\Omega)}^{\mathfrak{s}} \\
 & \hspace{2cm} + C \varepsilon^{-\gamma\mathfrak{s}-\frac{\mathfrak{s}}{2\mathfrak{q}}} \|\mathbf{M}_{\alpha}(|\mathbb{P}(x)\mathbf{G}|^{p})\|_{L^{\mathfrak{q},\mathfrak{s}}_{\mu}(\Omega)}^{\mathfrak{s}}.\label{LR-2}
\end{align}
To obtain~\eqref{est:B}, we simply fix $\varepsilon$ small enough in~\eqref{LR-2} such that $C \varepsilon^{\frac{\mathfrak{s}}{2\mathfrak{q}}}< \frac{1}{2}$. The same manner can be performed for the case $\mathfrak{s}=\infty$ to complete the proof.
\end{proof}

\begin{theorem}\label{theo:C}
Let $\mathbb{P}$ be a matrix weight satisfying~\eqref{cond-PI}; $\mathbf{F} \in L^{p}_{\omega}(\Omega)$ and $\mathsf{g} \in W^{1,p}_{\omega}(\overline{\Omega})$ for $p \in (1,\infty)$. Assume further that $u \in \mathsf{g} + W^{1,p}_{0,\omega}(\Omega)$ is a weak solution to~\eqref{Eq-0}. Given $\alpha \in [0,n)$, two weights $\mu \in \mathrm{A}_{\infty}$, $\nu \in L^1_{\mathrm{loc}}(\mathbb{R}^+;\mathbb{R}^+)$ and a function $\Sigma$ defined by~\eqref{def:Gam} satisfies~\eqref{cond-I} with two constants $c_1, c_2$. Then, for every $0<\mathfrak{q}<\infty$ and $0<\mathfrak{s}\le\infty$, there exists a constant $\kappa = \kappa(\mathtt{data},\mathfrak{q},\mathfrak{s},c_1,c_2)>0$ such that 
\begin{align}\label{est:theo-C}
\|\mathbf{M}_{\alpha}(|\mathbb{P}\nabla u|^{p})\|_{L^{\mathfrak{q},\mathfrak{s}}_{\mu,\nu}(\Omega)}  \le C\left\|\mathbf{M}_{\alpha}(|\mathbb{P}\mathbf{G}|^{p})\right\|_{L^{\mathfrak{q},\mathfrak{s}}_{\mu,\nu}(\Omega)},
\end{align}
 if $\log \mathbb{P}$ satisfies $(\kappa,r_0)$-small-$\log$-$\mathrm{BMO}$ condition and $\Omega$ is $(\kappa,r_0)$-Lipschitz domain for some $r_0>0$. Here, $C$ is a positive constant depending on $\mathfrak{q},\mathfrak{s},\mathtt{data}$.
\end{theorem}
\begin{proof}[Proof of Theorem~\ref{theo:C}]
Thanks to Theorem~\ref{theo:A}, for every $\theta>0$ and $\varepsilon>0$ small enough, one can find $\gamma = \gamma(\mathtt{data},\alpha)$ and $\kappa = \kappa(\mathtt{data},\theta,\alpha,\varepsilon)$ such that if $|\log \mathbb{P}|_{\mathrm{BMO}} \le \kappa$ and $\Omega$ is $(\kappa,r_0)$-Lipschitz for some $r_0>0$, then 
\begin{align}\notag
 & d^{\mu}_{\alpha}(|\mathbb{P}(x)\nabla u|^{p};\lambda) \le C \varepsilon d^{\mu}_{\alpha}(|\mathbb{P}(x)\nabla u|^{p};\varepsilon^{\theta}\lambda) + d^{\mu}_{\alpha}(|\mathbb{P}(x)\mathbf{G}|^{p};\varepsilon^{\theta+\gamma}\lambda)
\end{align}
for all $\lambda \in \mathbb{R}^+$. Thanks to Lemma~\ref{lem-Sig}, this inequality implies to
\begin{align}\label{ineq-d-1}
 \Sigma\big(d^{\mu}_{\alpha}(|\mathbb{P}(x)\nabla u|^{p};\lambda)\big) &\le c_2 \left[c_1 (C \varepsilon)^{\log_2 c_1} \Sigma \big(d^{\mu}_{\alpha}(|\mathbb{P}(x)\nabla u|^{p};\varepsilon^{\theta}\lambda)\big) \right. \notag \\
 & \qquad \qquad  \left. + \Sigma \big(d^{\mu}_{\alpha}(|\mathbb{P}(x)\mathbf{G}|^{p};\varepsilon^{\theta+\gamma}\lambda)\big) \right],
\end{align}
for $\varepsilon$ satisfying $C\varepsilon<1/2$. For simplicity, let us denote
\begin{align*}
\mathcal{G}(\lambda) := \Sigma\big(d^{\mu}_{\alpha}(|\mathbb{P}(x)\nabla u|^{p};\lambda)\big) \mbox{ and } \mathcal{F}(\lambda) := \Sigma \big(d^{\mu}_{\alpha}(|\mathbb{P}(x)\mathbf{G}|^{p};\lambda)\big).
\end{align*}
We can rewrite~\eqref{ineq-d-1} as below
\begin{align}\notag
\mathcal{G}(\lambda) &\le C \varepsilon^{\log_2 c_1} \mathcal{G}(\varepsilon^{\theta}\lambda) + C \mathcal{F}(\varepsilon^{\theta+\gamma}\lambda),
\end{align}
which yields to
\begin{align}
\left(\int_{\mathbb{R}^+} \mathfrak{q} \big[\lambda^{\mathfrak{q}}\mathcal{G}(\lambda)\big]^{\frac{\mathfrak{s}}{\mathfrak{q}}} \frac{d\lambda}{\lambda}\right)^{\frac{1}{\mathfrak{s}}} &\le C \varepsilon^{\frac{1}{\mathfrak{q}}\log_2 c_1} \left(\int_{\mathbb{R}^+} \mathfrak{q} \big[\lambda^{\mathfrak{q}}\mathcal{G}(\varepsilon^{\theta}\lambda)\big]^{\frac{\mathfrak{s}}{\mathfrak{q}}} \frac{d\lambda}{\lambda}\right)^{\frac{1}{\mathfrak{s}}} \notag \\
& \qquad \qquad + C \left(\int_{\mathbb{R}^+} \mathfrak{q} \big[\lambda^{\mathfrak{q}}\mathcal{F}(\varepsilon^{\theta+\gamma}\lambda)\big]^{\frac{\mathfrak{s}}{\mathfrak{q}}} \frac{d\lambda}{\lambda}\right)^{\frac{1}{\mathfrak{s}}}. \label{ineq-d-1b} 
\end{align}
By changing of variables and using notion of quasi-norm in~\eqref{GL-norm}, we obtain from~\eqref{ineq-d-1b} that
\begin{align}\notag
 \|\mathbf{M}_{\alpha}(|\mathbb{P}(x)\nabla u|^{p})\|_{L^{\mathfrak{q},\mathfrak{s}}_{\mu,\nu}(\Omega)} & \le C \varepsilon^{\frac{1}{\mathfrak{q}}\log_2 c_1 - \theta}  \|\mathbf{M}_{\alpha}(|\mathbb{P}(x)\nabla u|^{p})\|_{L^{\mathfrak{q},\mathfrak{s}}_{\mu,\nu}(\Omega)} \\
 & \hspace{2cm} + C \varepsilon^{-\theta-\gamma} \|\mathbf{M}_{\alpha}(|\mathbb{P}(x)\mathbf{G}|^{p})\|_{L^{\mathfrak{q},\mathfrak{s}}_{\mu,\nu}(\Omega)}.\label{ineq-d-2}
\end{align}
We also obtain~\eqref{ineq-d-2} by similar ways for the remain case $\mathfrak{s} = \infty$. To finish the proof, we just choose $\theta<\frac{1}{\mathfrak{q}}\log_2 c_1$ and fix $\varepsilon$ small enough such that $C \varepsilon^{\frac{1}{\mathfrak{q}}\log_2 c_1 - \theta}<1/2$. We then obtain~\eqref{est:theo-C}.
\end{proof}

\begin{theorem}\label{theo:D}
Let $\alpha \in [0,n)$; $\mathbb{P}$ be a matrix weight satisfying~\eqref{cond-PI}; $\mathbf{F} \in L^{p}_{\omega}(\Omega)$ and $\mathsf{g} \in W^{1,p}_{\omega}(\overline{\Omega})$ for $p \in (1,\infty)$. Assume further that $u \in \mathsf{g} + W^{1,p}_{0,\omega}(\Omega)$ is a weak solution to~\eqref{Eq-0}. Assume further that $\upsilon \in (0,n)$ and $\psi: \Omega \times \mathbb{R}^+ \to \mathbb{R}^+$ satisfies the following condition 
\begin{align}\label{cond-psi-2}
\psi(x,2t) \le 2^{\upsilon} \psi(x,t), \quad \mbox{ for all } x \in \Omega \mbox{ and } t > 0.
\end{align}
Then, for every $\mathfrak{q}>0$, there exists a constant $\kappa = \kappa(\mathtt{data},\mathfrak{q},\upsilon)>0$ such that
\begin{align}\label{est:theo-D}
\|\mathbf{M}_{\alpha}(|\mathbb{P}\nabla u|^{p})\|_{\mathrm{M}^{\mathfrak{q},\psi}(\Omega)}  \le C \left\|\mathbf{M}_{\alpha}(|\mathbb{P}\mathbf{G}|^{p})\right\|_{\mathrm{M}^{\mathfrak{q},\psi}(\Omega)}.
\end{align}
 if $\log \mathbb{P}$ satisfies $(\kappa,r_0)$-small-$\log$-$\mathrm{BMO}$ condition and $\Omega$ is $(\kappa,r_0)$-Lipschitz domain for some $r_0>0$. Here, $C$ is a positive constant depending on $\mathfrak{q},\upsilon,\mathtt{data}$.
\end{theorem}
\begin{proof}[Proof of Theorem~\ref{theo:D}]
Let $y \in \Omega$ and $0<\varrho<D_0$, we consider $\mathbb{K}$ as
\begin{align*}
\mathbb{K} & := \frac{1}{\psi(y,\varrho)}\int_{\Omega(y,\varrho)}[\mathbf{M}_{\alpha}(|\mathbb{P}(x)\nabla u|^{p})]^{\mathfrak{q}}dx.
\end{align*}
Since $\chi_{B(y,\varrho)}(x) \le \left[\mathbf{M} \chi_{B(y,\varrho)}\right](x) \le \left[\mathbf{M} \chi_{B(y,\varrho)}\right]^{\sigma}(x) \le 1$ for $x \in \mathbb{R}^n$ a.e. and for all $\sigma \in (0,1)$, then it follows that
\begin{align}\notag 
\mathbb{K}  & \le \frac{1}{\psi(y,\varrho)} \int_{\mathbb{R}^n}[\chi_{\Omega}(x)\mathbf{M}_{\alpha}(|\mathbb{P}(x)\nabla u|^{p})]^{\mathfrak{q}} \left[\mathbf{M} \chi_{B(y,\varrho)}\right]^{\sigma}(x) dx.
\end{align}
Thanks to~\cite[Proposition 2]{CR80}, we conclude that $\left[\mathbf{M} \chi_{B(y,\varrho)}\right]^{\sigma} \in \mathrm{A}_{\infty}$. Hence, it is possible to apply Theorem~\ref{theo:A} with $\mathfrak{s}=\mathfrak{q}$ and $\mu = \left[\mathbf{M} \chi_{B(y,\varrho)}\right]^{\sigma}$. Then, we get that
\begin{align}\label{est-D-2}
\mathbb{K} & \le \frac{C}{\psi(y,\varrho)} \int_{\mathbb{R}^n} \left[\chi_{\Omega}(x)\mathbf{M}_{\alpha}(|\mathbb{P}(x)\mathbf{G}|^{p})\right]^{\mathfrak{q}} \left[\mathbf{M} \chi_{B(y,\varrho)}\right]^{\sigma}(x) dx.
\end{align}
Using the notation in~\eqref{Omega-j}, we have the following decomposition
\begin{align*}
\mathbb{R}^n = B(y,2\varrho) \cup \left(\bigcup_{j=1}^{\infty} \Omega_{j}^{\varrho}(y) \right).
\end{align*}
We may rewrite~\eqref{est-D-2} as follows
\begin{align} 
\mathbb{K} & \le \frac{C}{\psi(y,\varrho)}\int_{B(y,2\varrho)} \left[\chi_{\Omega}(x)\mathbf{M}_{\alpha}(|\mathbb{P}(x)\mathbf{G}|^{p})\right]^{\mathfrak{q}} \left[\mathbf{M} \chi_{B(y,\varrho)}\right]^{\sigma}(x) dx \notag \\
& \qquad + \frac{C}{\psi(y,\varrho)} \sum_{j=1}^{\infty}\int_{\Omega_{j}^{\varrho}(y)} \left[\chi_{\Omega}(x)\mathbf{M}_{\alpha}(|\mathbb{P}(x)\mathbf{G}|^{p})\right]^{\mathfrak{q}} \left[\mathbf{M} \chi_{B(y,\varrho)}\right]^{\sigma}(x) dx.\notag
\end{align}
Applying condition~\eqref{cond-psi-2} and inequality~\eqref{gmu} in Lemma~\ref{lem:Mor}, one obtains that
\begin{align}\label{est-D-4}
\mathbb{K} & \le C \left(\frac{2^{\upsilon}}{\psi(y,2\varrho)}\int_{B(y,2\varrho)} \left[\chi_{\Omega}(x)\mathbf{M}_{\alpha}(|\mathbb{P}(x)\mathbf{G}|^{p})\right]^{\mathfrak{q}} dx\right) \notag \\
& \qquad +  C \sum_{j=1}^{\infty} 2^{-(j-1)n\sigma} \left(\frac{2^{(j+1)\upsilon}}{\psi(y,2^{j+1}\varrho)} \int_{B(y,2^{j+1}\varrho)} \left[\chi_{\Omega}(x)\mathbf{M}_{\alpha}(|\mathbb{P}(x)\mathbf{G}|^{p})\right]^{\mathfrak{q}} dx\right) \notag \\
& \le 2^{\upsilon+n\sigma} C \left[1 + \sum_{j=1}^{\infty} \left(2^{\upsilon-n\sigma}\right)^j\right] \|\mathbf{M}_{\alpha}(|\mathbb{P}(x)\mathbf{G}|^{p})\|_{\mathrm{M}^{\mathfrak{q},\psi}(\Omega)}^{\mathfrak{q}}.
\end{align}
Since $\upsilon < n$, it is possible to fix $\sigma \in (0,1)$ in~\eqref{est-D-4} such that
$2^{\upsilon-n\sigma} < 1$. It implies that the series $\sum_{j=1}^{\infty} \left(2^{\upsilon-n\sigma}\right)^j$ is finite. The proof of~\eqref{est:theo-D} is now complete.
\end{proof}

\section*{Acknowledgement}
This research is funded by Vietnam National Foundation for Science and Technology Development (NAFOSTED), Grant Number: 101.02-2021.17.

\section*{Conflict of Interest}
The authors declared that they have no conflict of interest.

\section*{Declarations}
Data sharing not applicable to this article as no datasets were generated or analysed during the current study.\\

\end{document}